\documentclass{article}%
\usepackage{graphicx}
\usepackage{amsmath}%
\usepackage{amsfonts}%
\usepackage{amssymb}

\usepackage{amssymb}
\usepackage{amsbsy}
\newcommand{\no}{\noindent}
\newcommand{\be}{\begin{equation}}
\newcommand{\ee}{\end{equation}}
\newtheorem{prop}{Proposition}[section]

\def\Z{\mathbb Z}

\def\Q{\mathbb Q}
\def\H{\mathcal H}

\def\R{\mathbb R}
\def\C{\mathbb C}

\def\F{\mathcal F}

\title{Twisted K-theory and cohomology}
\author{Michael Atiyah\\ \small School of Mathematics, University of Edinburgh,\\ \small James Clerk Maxwell Building,\\ \small King's Buildings, Edinburgh EH9 3JZ, U.K. \and Graeme Segal\\ \small All Souls College, Oxford OX1 4AL, U.K.}
\begin{document}
\maketitle
\section{\bigskip Introduction}

\ \ \ \  In a previous paper [AS] we introduced the twisted $K$-theory of a space $X$, with the twisting corresponding to a 3-dimensional cohomology class of $X.$ \ If the class is
of finite order then the theory, particularly when tensored with the rationals $\Q,$ is very
similar to ordinary $K$-theory. \ For elements of infinite order, however, 
there are substantial differences. \ In particular the relations with
cohomology theory are very different, and it is these that are the
subject of this paper.\bigskip

 There are three different ways in which ordinary $K$-theory and
cohomology are related:

\begin{enumerate}
\item[(1)] by the Atiyah-Hirzebruch spectral sequence,

\item[(2)] by the Chern character,

\item[(3)] by the Chern classes.
\end{enumerate}

\noindent In the first and third of these integral cohomology is involved,
while in the second only rational cohomology enters.\bigskip

 In the twisted case we shall examine the counterparts of all three
of these relations, finding essential new features in all of them. \ We shall
also briefly examine {\it operations} in twisted $K$-theory, the analogue of the
well-known $\lambda$-ring operations in ordinary $K$-theory, because these have some bearing
on the cohomology questions.\bigskip

 Twisted $K$-theory arises because $K$-theory admits natural automorphisms
given by tensoring with line-bundles. \ We review this briefly in Section 2,
and then use the basic formulae in Sections 3 and 4 to compute the relevant
differentials of the appropriate spectral sequences. \ The first key result is
the explicit identification of the differential $d_{3}$ of the
Atiyah-Hirzebruch spectral sequence already forecast in [AS] (cf. also [R]). \ In Section
5 we enter new territory by showing that \emph{the higher differentials}
$d_{5},d_{7},...$ are \emph{non-zero even over the rationals or reals}, and
given by Massey products. \ (In the Appendix we give explicit examples, pointed out to us by Elmer Rees, of spaces for which the higher differentials do not vanish.) \ In Section 6 we introduce (over the reals) the
appropriate twisted cohomology. \ A key feature here is that the theory is not
integer graded but (like $K$-theory) is filtered (and only $\operatorname{mod}%
2$ graded). \ In terms of differential forms such a theory has already
appeared, at least implicitly, in the physics literature. \ This twisted
cohomology has its own \textquotedblleft Atiyah-Hirzebruch spectral
sequence\textquotedblright. \ Comparing this with twisted $K$-theory we show
in Section 7 the existence of a \emph{Chern character} in the twisted
context.\bigskip

 The last three sections are then devoted to Chern classes for
twisted $K$-theory. \ First we show that the ring of characteristic classes
(in ordinary real cohomology) is the subring of the algebra of Chern classes
invariant under tensoring by line-bundles (Proposition (8.8)). \ We then
examine how these might be lifted to twisted characteristic classes.\bigskip

 In Section 9 we discuss a very geometric set of characteristic
classes studied by Koschorke [K], and show how they fit into the twisted
theory.\bigskip

 Finally in Section 10 we discuss internal operations in twisted
$K$-theory.\bigskip

 These last sections identify a number of purely algebraic problems
which arise, and, though it would be nice to know the answers, we have not felt
that there was sufficient motivation or application as yet to justify
embarking on the possibly lengthy algebraic investigation.\bigskip

 In a differential-geometric context there is a very rich theory
relating curvature to differential forms representing the (real) Chern
classes. \ These are of course of special interest to both geometers and
physicists, and they can lead to more precise local formulae. It is natural to ask for a similar theory for twisted $K$-theory,
but we have given only a brief sketch of it in Section 7. \ Fuller treatments of this subject can be found in [BCMMS] and in unpublished work of Bunke. \  (Cf. also [TX].) \ There are in fact
some significant technical difficulties, at least from our approach.
\ Although the space of Fredholm operators is an elegant model for $K$-theory
from the point of view of topology, it is too large or too crude for
differential geometry. \ This point arises clearly in Connes's non-commutative
geometry, where $C^{\ast}$-algebras are too crude to provide explicit cocycle
formulae, and the theory has to be refined to incorporate some differential
calculus.  \ It seems to us that Connes's theory should provide the right
framework to deal with curvature and Chern forms.\bigskip

\section{\bigskip The action of the automorphism group}

\ \ \ \  In our previous paper [AS] \ we explained briefly how automorphisms
of a cohomology theory lead to twisted versions of that theory. \ In
particular, for $K$-theory, we studied the automorphisms given by tensoring
with line bundles, and the twisting that this determines. \ In this section we
shall examine the cohomological implications of these automorphisms.\bigskip

 For a compact space $X$ we have the group Pic$\left(  X\right)  $
formed by (classes of) line bundles with the tensor product as multiplication.
\ The first Chern class gives an isomorphism 
\begin{equation}
c_{1}:\text{Pic}\left(  X\right)  \longrightarrow H^{2}\left(  X;Z\right)
\tag{2.1}%
\end{equation}
Moreover Pic$\left(  X\right)  $ is naturally a subgroup of the multiplicative
group of the ring $K\left(  X\right)  ,$ and hence it acts by group
homomorphisms on $K\left(  X\right)  .$ \ This action is functorial in $X$, and
so it can equally well be described in terms of classifying spaces.

\bigskip

 The classifying space of $K$-theory can be taken to be Fred$\left(
\H\right)  ,$ the space of Fredholm operators on Hilbert space $\H,$ while the
classifying space of Pic can be taken to $PU\left(  \H\right)  ,$ the
projective unitary group of $\H.$ \ The group structure of $PU\left(  \H\right)
$ gives the group structure of Pic$\left(  X\right)  $.

\begin{prop} The action of  \ ${\rm Pic}(X)$ on
$K\left(  X\right)  $ is induced by the conjugation action of $PU\left(
\H\right)  $ on ${\rm Fred}\left(  \H\right)  .$
\end{prop}

\no{\bf Proof}\ \ If $L$ is a line bundle on $X$, and an element $\xi$ of $K(X)$ is represented by a family $F=\{F_x\}$ of Fredholm operators in a fixed Hilbert space $\H$, then $L\otimes \xi \in K(X)$ is represented by the family of Fredholm operators $1\otimes F=\{1\otimes F_x\}$ in the Hilbert bundle $L\otimes \H$ on $X$. This bundle can be trivialized by an isomorphism $u:L\otimes \H\to X \times \H$, which exists and is unique up to homotopy because $U(\H)$ is contractible. So the class $L\otimes \xi$ can be represented by the family $u\circ(1\otimes F)\circ u^{-1}$ in the fixed space $\H$. This proves the proposition, for we can regard $u$ as a map $u:X\to PU(\H)$, and the line bundle on $X$ pulled back from $PU(\H)$ by $u$ is precisely $L$. $\Box$ 

\bigskip

 We will now consider the effect on cohomology of the map
\begin{equation}
PU\left(  \H\right)  \times\text{Fred}\left(  \H\right)  \longrightarrow
\text{Fred}\left(  \H\right) \tag{2.2}%
\end{equation}
or equivalently the effect on cohomology of the multiplication
\begin{equation}
\text{Pic}\left(  X\right)  \times K\left(  X\right)  \longrightarrow K\left(
X\right) \tag{2.3}%
\end{equation}

\noindent In other words we need to describe the Chern classes of a
vector bundle $L\otimes E$ in terms of the Chern classes of the vector bundle
$E$ and the line bundle $L.$ \ These are given by universal polynomials, and
since the classifying spaces in (2.2) have no torsion we lose nothing by working with
rational cohomology.

\bigskip

 As usual we define the Chern classes $c_{n}\left(  E\right)  $ of an $N$-dimensional bundle $E$ in
terms of variables $x_{j}$ by
\[
\sum\limits_{n=0}^{N}c_{n}=\prod\limits_{j=1}^{N}\left(  1+x_{j}\right),
\]
and the power sums $s_{n}$ by
\[
s_{n}=\sum\limits_{j=1}^{N}x_{j}^n.
\]
In particular, $s_0(E)=\sum x^0_j$ is the {\it dimension} of $E$. \ As a function on Fred$(\H)$ it becomes the {\it index} map, labelling the connected components of Fred$(\H)$.

\bigskip

 The Chern character is given by
\[
{\rm ch}\;E=\sum\limits_{n=0}^{\infty}\dfrac{s_{n}}{n!},%
\]
and it extends to a ring homomorphism
\[
{\rm ch}:K^{\ast}\left(  X\right)  \longrightarrow H^{\ast}\left(  X,\Q\right)
\]
which becomes an isomorphism after tensoring with $\Q$.

\bigskip

 If $L$ is a line-bundle with $c_{1}\left(  L\right)  =u,$ and $E$ is
a vector bundle, then the Chern character of $L\otimes E$ is
given by
\[
{\rm ch}\left(  L\otimes E\right)  =e^{u}\;{\rm ch}E.
\]
This implies the same formula for any element $E$ in $K\left(  X\right)  .$
\ If $s_{n}$ and $s_{n}\left(  u\right)  $ denote the power sums for $E$ and
$L\otimes E$ then
\begin{equation}
\sum\dfrac{s_{n}\left(  u\right)  }{n!}=\left(  1+u+\dfrac{u^{2}}%
{2!}+...\right)  \sum\dfrac{s_{n}}{n!}\tag{2.4}%
\end{equation}
so that
\begin{equation}
s_{n}\left(  u\right)  =s_{n}+ns_{n-1}u+\left(
\genfrac{}{}{0pt}{0}{n}{2}%
\right)  s_{n-2}u^{2}+...+s_{0}u^{n}.\tag{2.5}%
\end{equation}
This formula can be interpreted as giving the homomorphism in cohomology
induced by the map (2.2). 
\bigskip

 In the next section we shall deduce various consequences from the
formula (2.5).

\section{The universal fibration}

\ \ \ \  The group action (2.2) of $PU\left(  \H\right)  $ on Fred$(\H)$ gives
rise to a \textquotedblleft universal bundle\textquotedblright\ with fibre
Fred$\left( \H\right)  $ over the classifying space $BPU\left(  \H\right)  .$
\ This classifies $PU\left(  \H\right)  $-bundles over a space $X$ as homotopy classes of maps
$X\rightarrow BPU\left(  \H\right)  .$ \ Note that a $PU\left(  \H\right)
$-bundle determines a $P(\H)$-bundle and conversely; in [AS] we worked
primarily with $P\left(  \H\right)  $-bundles.

\bigskip

 Since $PU\left(  \H\right)  =U\left(  \H\right)  /U\left(  1\right)  $,
it is an Eilenberg-MacLane space $K\left(  \Z,2\right)  $, and so its
classifying space $BPU\left(  \H\right)  $ is a $K\left(  \Z,3\right)  .$
\ Thus, as observed in [AS] Prop.2.1, isomorphism classes of $PU\left(  \H\right)
$-bundles over $X$ are classified by elements of $H^{3}\left(  X,\Z\right)  .$

\bigskip

 In the notation of [AS] our universal bundle over $BPU\left(
\H\right)  $ with fibre Fred$\left(  \H\right)  $ is written as Fred$\left(
P\right)  $, where $P$ is the universal $P\left(  \H\right)$ --- or 
$PU\left(  \H\right)$ --- bundle over $BPU\left(  \H\right)  .$
\ We now want to study the rational cohomology of this fibration, using the
formula (2.5). \ Note that because the index is preserved by the group action the connected components
Fred$_k\left(  P\right)  $ of Fred$\left(  P\right)  $ also correspond to the index.

\bigskip

 Now, over the rationals, $K\left(  \Z,3\right)  $ is indistinguishable from a
3-sphere. \ To have an explicit map $S^3\to BPU(\H)$, note that, for any group $G,$ we have an
inclusion of its suspension $\Sigma G$ into $BG$ (corresponding to $G$-bundles with just one transition
function, when $\Sigma G$ is regarded as the union of two cones on $G$). \ Hence we have inclusions
\begin{equation}%
\Sigma PU(\H) \to BPU(\H)
\tag{3.1}%
\end{equation}%
\begin{equation}
S^{3}=\Sigma (\C P^1) \to \Sigma(\C P^{\infty}) \simeq \Sigma (PU(\H)). \tag{3.2}%
\end{equation}
The composite map
\begin{equation}
S^{3}\longrightarrow BPU\left(  \H\right) \tag{3.3}%
\end{equation}
induces an isomorphism on $H^{3}\left(\   \;;\Z\right)  ,$ while all higher
cohomology groups of $BPU\left(  \H\right)  \simeq K\left(  \Z,3\right)  $ are finite.

\bigskip

 Consequently, if  $Y_{k}$ denotes the restriction of the bundle
Fred$_{k}\left(  P\right)  $ to $S^{3},$ our task is to examine the Leray
spectral sequence (over $\Q)$ of $Y_{k}$. \ Since the base is $S^{3}$ the only
non-zero differential of the spectral sequence is $d_{3},$ and the spectral
sequence  reduces to the Wang exact sequence which is valid for any
fibration over a sphere $S^{m},$ with  total space $Y.$ \ This is
just the exact sequence of the pair $\left(  Y,F_{0}\right)  ,$ where $F_0$ is the fibre over the base-point $0\in
S^{m}$.
\begin{equation}
...\longrightarrow H^{q-1}\left(  F_{0}\right)  \overset{\delta}%
{\longrightarrow}H^{q}\left(  Y,F_{0}\right)  \longrightarrow H^{q}\left(
Y\right)  \longrightarrow H^{q}\left(  F_{0}\right)  \longrightarrow
...\tag{3.4}%
\end{equation}
If $\infty\in S^{m}$ is the antipodal base-point, and $F_{\infty}$ is the fibre over it, then we have a suspension
isomorphism
\[
\sigma:H^{q-m}\left(  F_{\infty}\right)  \longrightarrow H^{q}\left(
Y,F_{0}\right)  .
\]
Substituting this in (3.4) we get the Wang exact sequence
\begin{equation}
...\longrightarrow H^{q-1}\left(  F_{0}\right)  \overset{d}{\longrightarrow
}H^{q-m}\left(  F_{\infty}\right)  \longrightarrow H^{q}\left(  Y\right)
\longrightarrow H^{q}\left(  F_{0}\right)  \longrightarrow...\; ,\tag{3.5}%
\end{equation}
where (identifying $H^m(S^m,H^{q-m}(F_{\infty}))$ with $H^{q-m}(F_{\infty})$)
\begin{equation}
d=\sigma^{-1}\delta\tag{3.6}%
\end{equation}
is  the differential $d_{m}$ of the Leray spectral sequence
\begin{gather*}
d_{m}:E_{2}^{0,q-1}\longrightarrow E_{2}^{m,q-m}.\\
\end{gather*}

 Because the base of the fibration is simply connected we can identify $H^*(F_{\infty})$ canonically with $H^*(F_0)$, and then the homomorphism $d:H^*(F_0)\to H^*(F_0)$ is a {\it derivation}. \ We can explicitly calculate $d$ from the homomorphism $\varphi
^{\ast}$ induced by the gluing map
\begin{equation}
\phi:S^{m-1}\times F_{\infty}\longrightarrow S^{m-1}\times F_{0}\tag{3.7}%
\end{equation}
which (homotopically) defines the fibration $Y\rightarrow X.$ \ Using (3.6)
one finds the formula 
\begin{equation}
\phi^{\ast}\left(  1\otimes x\right)  = (1\otimes x) +\left(  s\otimes dx\right), \tag{3.8}%
\end{equation}
where $x\in H^{q-1}\left(  F_{0}\right)$ and 1 and $s$ are the generators of
$H^{0}\left(  S^{m-1}\right)  $ and $H^{m-1}\left(  S^{m-1}\right)$ respectively. \ (Notice that

 $d$ is a derivation \ \ $\Longleftrightarrow$ \ \  $\varphi ^*$ is a ring homomorphism.)

\bigskip

 We now take $m=3$ and $F_0=$\ Fred$_{k}\left(  \H\right)  ,$ with the
fibration over $S^{3}$ induced by the inclusion of $S^{3}$ in $BPU\left(
\H\right)  $ as in (3.3). \ The map $\phi$ is then the restriction of the
multiplication map (2.2) to
\[
S^{2}=\C P^{1}\subset \C P^{\infty}\simeq PU\left(  \H\right)  .
\]
Cohomologically this means that, in formula (2.5), we drop all powers of $u$
except the linear term. \ In particular this implies that only for $n=1$ does
the index $k=s_{0}$ of the component Fred$_{k}\left(  \H\right)  $ enter:
\[
s_{1}\left(  u\right)  =s_{1}+ku.
\]
Putting all this together, we see that we have established:\bigskip

\noindent\textbf{Proposition (3.9) \ }\emph{\ In the spectral sequence for the
universal fibration over }$S^{3},$\emph{\ with fibre }Fred$_{k}\left(
\H\right)  ,$\emph{\ the differential }$d_{3}$\emph{\ is given on the
cohomology generators }$s_{n}$\emph{\ by the formulae:}
\begin{align*}
d_{3}s_{n}  & =ns_{n-1}w\;\text{for\ \ }n\geq2\\
d_{3}s_{1}  & =s_0w=kw,
\end{align*}
\emph{where }$w$ \emph{is the generator of }$H^{3}\left(  S^{3};\Z\right)  .$

\bigskip

 The coefficients in (3.9) can be taken to be the integers. \ Over the
rationals, as noted before, $S^{3}$ is equivalent to $BPU(\H)$ and so the
formula in (3.9) will hold over $\Q$ for all fibrations. \ If the fibration
over $X$ is induced by a class $\eta\in H^{3}\left(  X;\Z\right)  ,$ then the spectral sequence for
the induced fibration over $X,$ has differential $d_{3}
$ given by
\begin{align*}
d_{3}s_{n}  & =ns_{n-1}\eta\;\;\;\;\;n\geq2\\
d_{3}s_{1}  & =k\;\eta
\end{align*}
In particular,  the spectral sequence yields the exact sequence
\[
H^{2}\left(  \text{Fred}_{k}\right)  \overset{d_{3}}{\longrightarrow}%
H^{3}\left(  X\right)  \overset{\pi^{\ast}}{\longrightarrow}H^{3}\left(
\text{Fred}_{k}\left(  P\right)  \right),
\]
where $\pi$ is the fibre projection. Hence $\pi^{\ast}d_{3}=0$ \ and so  $\pi^{\ast}\eta=0$ if $k \neq 0$. But if the bundle possesses a section then $\pi^*$ is injective, so if $\eta \neq 0 $ rationally --- i.e. if $\eta \in H^3(X;\Z)$ has infinite order --- we must have $k=0$.  \ This shows that, as forecast in [AS]\ \S 3 ~Remark~(iv)],
\emph{for twistings\ of infinite order there are no sections with non-zero
index. \ }In other words, the index map $K_{P}\left(  X\right)
\rightarrow \Z$ is zero for $[P]$ of infinite order.

\bigskip

 In  subsequent sections we shall pursue further the cohomological
implications of the formula in (3.9), but since our main concern is to relate
twisted $K$-theory to cohomology we next turn to this.

\section{The Atiyah-Hirzebruch spectral sequence}

\ \ \ \  For ordinary $K$-theory there is the Atiyah-Hirzebruch spectral
sequence [AH] which relates it to cohomology. \ As already briefly explained in [AS]  Prop (4.1), there is an analogous spectral sequence in the
twisted case. \ This has the same $E_{2}$ term
\[
E_{2}^{p,q}=H^{p}\left(  X,K^{q}\text{(point)}\right)
\]
as in the untwisted case, but converges to the graded group $E_{\infty}^{p,q}$ associated to a
filtration on $K_{P}^{\ast}\left(  X\right)  .$ \ Because $K^q({\rm point}) =0$ when $q$ is odd the even differentials $d_{2},d_{4},...$ are zero. \ In the
untwisted case the first non-zero differential $d_{3}$ is (see [AH]) the Steenrod
operation
\[
Sq^{3}_{\Z}:H^{p}\left(  X;\Z\right)  \longrightarrow H^{p+3}\left(  X;\Z\right)  .
\]
Our first task is to find $d_{3}$ in the twisted case. \ We shall show
that
\begin{equation}
d_{3}\left(  x\right)  =Sq^{3}_{\Z}\left(  x\right)  -\eta x_{{}}\tag{4.1}%
\end{equation}
where $\eta=\eta\left[  P\right]  $ is the class of $P.$

\bigskip

 By a standard argument in homotopy theory, the only universal operations on $H^p(\;\; ;\Z)$,
raising dimension by 3, defined for spaces with a given class $\eta\in
H^{3},$ are given by elements of%
\[
H^{p+3}\left(  K\left(  \Z,p\right)  \times K\left(  \Z,3\right)  ;\Z\right)
\]
This group is isomorphic to

\[H^{p+3}(K(\Z,p))\oplus H^{p+3}(K(\Z,3)) \oplus \Z,\]
where the third summand is generated by the product of the generators of $H^p(K(\Z,p))$ and $H^3(K(\Z,3))$. We deduce that the operation
$d_{3}$ must be of the form%
\begin{equation}
d_{3}\left(  x\right)  =Sq^{3}_{\Z}\left(  x\right)  +b\eta x_{{}}\tag{4.2}%
\end{equation}
for some $b\in \Z$ , for the formula must agree with the
untwisted case when $\eta=0$, and the operation must vanish when $x=0$. It remains to compute the
integer $b.$ \ For this it is sufficient to consider the spectral sequence for the twisted $K$-theory of $S^3$, with the generator $\eta$ of $H^3(S^3;\Z)$ as the class of the twisting. \ We
just have to compute
\begin{equation}
d_{3}:H^{0}\left(  S^{3}\right)  \longrightarrow H^{3}\left(  S^{3}\right)
\tag{4.3}%
\end{equation}
and show that $d_{3}\left(  1\right)=-\eta  $. \  This proceeds  analogously to the calculation of the
$d_{3}$ of Section 3.
\bigskip 

The spectral sequence for a space $X$ is obtained by applying the functor $K^*_P$ to the filtered space $X$, filtered by its skeletons. \ When $X=S^3$ the filtration has just a single step $X^0 = X^1 = X^2 \subset X^3 =X$, where $X^0$ is the basepoint, and so the spectral sequence reduces to the long exact sequence for the pair $(X,X^0)$.\  The differential (4.3) is easily seen to be the boundary map\  $K^0_P(X^0)\to K^1_P(X,X^0)$.\ Equivalently, it is the boundary map\  $K^0_P(D_0)\to K^1_P(S^3,D_0)$, when $S^3$ is written as the union of two hemispheres $D_0$ and $D_{\infty}$ intersecting in the equatorial sphere $S^2$. \ This map fits into the commutative square
\[\begin{array}{rcc}K^0_P(D_0)\ &\to \ &K^1_P(S^3,D_0)\\
 \downarrow \ && \ \downarrow\\
K^0_P(S^2) \ &\to \ &K^1_P(D_{\infty},S^2).
\end{array}\]
The differential of (4.3) is the passage from top left to bottom right in this diagram, when the groups are identified with $K^0(D_0)\cong \Z$ and $K^1(D_{\infty},S^2)\cong \Z$ by trivializing $P$ over $D_0$ and $D_{\infty}$ respectively. \ Now $K^0(S^2)\cong \Z \oplus \Z$, with generators 1 and $[L]$, where $L$ is the standard line bundle whose first Chern class is the generator of $H^2(S^2)$. \ The lower horizontal map takes 1 to 0 and $[L]$ to 1 when the trivialization from $D_{\infty}$ is used for both groups. \ The left-hand vertical map takes 1 to 1 if we use the $D_0$-trivialization on $S^2$, but takes 1 to $[L^{-1}]=2-[L]$ when $P$ is trivialised from $D_{\infty}.$ \ So, finally, the differential takes 1 to -1, as we claimed. \ We have thus established the result anticipated in [AS]
\bigskip

\noindent\textbf{Proposition (4.6) \ }\emph{In the Atiyah-Hirzebruch spectral
sequence for the functor }$K_{P}$\emph{\ the differential }$d_{3}$\emph{\ is
given by }
\[
d_{3}\left(  x\right)  =Sq^{3}_{\Z}\left(  x\right)  -\eta x,
\]
\emph{where }$\eta$\emph{\ is the class of }$P$\emph{\ in }$H^{3}\left(
X;\Z\right)  $

\section{The higher differentials}

\ \ \ \  For ordinary K-theory all the differentials of the Atiyah-Hirzebruch
spectral sequence are of finite order and so disappear on tensoring with the
rationals, showing that the graded group associated to a filtration of $K^{\ast}\otimes \Q$ is isomorphic to
$H^{\ast}\otimes \Q.$ \ If the twisting class $\eta\left(  P\right)  $ is of
finite order then the same  holds for $K_{P}^{\ast}$ , but if $\eta(P)$
is of infinite order then already, as shown in Section 3, $d_{3}$ can be of
infinite order. \ At first sight one might be tempted to think that this might
be all that happens, and  might expect
\begin{equation}
E_{4}\otimes \Q\cong E_{\infty}\otimes \Q.\tag{5.1}%
\end{equation}
However, as pointed out to us by Michael Hopkins, this is \textbf{not} the
case in general, although there are important cases where it does hold. \ The
purpose of this section is to explain precisely what does happen to the
spectral sequence when we tensor with $\Q.$

\bigskip

 As explained in Section 4, general homotopy theory  implies
that,\ in our category of spaces having a preferred class $\eta$ in $H^{3}$,
the only \emph{primary} cohomology operations over $\Q$ are given by products
with $\eta.$ \ Since $\eta$ is odd-dimensional, $\eta^{2}=0$ and so we only
get multiplication by $\eta,$ as in the $d_{3}$ for $K_{P}.$ \ But now let us
consider possible \emph{secondary }operations, i.e. ones defined not on $E_{2},$ but
on $E_{4}.$ \ To get $E_{4}$ we start from $E_{2}=E_{3}=H^{\ast}\left(
X,\Q\right)  $ and form first the $d_{3}$-cycles
\[
Z_{4}=\left\{  x\in H^{\ast}\left(  X\right)  \;\text{with\ }\eta x=0\right\}
\]
and then the boundaries
\[
B_{4}=\left\{  y\in H^{\ast}\left(  X\right)  \;\text{of the form\ }y=\eta
x\right\}
\]
giving
\[
E_{4}=Z_{4}/B_{4}.
\]
Now on $E_{4}=E_{5}$ there is a new universal operation that can be defined,
namely the \emph{Massey triple product}. \ Let us rapidly review the theory of
such products. \ We need to introduce cochains, and although it is
technically possible [Su] to do this over $\Q$ , \ it is simpler to assume our
spaces are manifolds and to work over $\R,$ where we can use differential forms
as cochains. \  In fact little is lost by this assumption, and 
for any application in physics the restriction is quite natural.

\bigskip

 In general, let $x,y,z$ be closed differential forms on $X$ of
degrees $p,q,r$ with
\[
\left[  x\right]  \left[  y\right]  =0\; ,\;\;\left[  y\right]  \left[
z\right]  =0,
\]
where [ \ \ ] denotes the de Rham cohomology  class. \ Thus
\[
xy=du\; ,\;\;yz=dv.
\]
Then the $\left(  p+q+r-1\right)  $ - form
\[
w=uz+\left(  -1\right)  ^{p-1}xv
\]
satisfies
\[
dw=duz-xdv=xyz-xyz=0.
\]
The cohomology class
\[
\left[  w\right]  =\left\{  x,y,z\right\}
\]
is called the {\it Massey triple product}. Because of the
choices of $u$ and $v$ , it is ambiguous up to the addition of arbitrary multiples of $x$ and $z$.

\bigskip

 There is a whole hierarchy of higher Massey products, and they play
a basic role in rational homotopy theory [Su]. \ Essentially they give all the
additional information (beyond the cohomology ring structure) that is needed to
determine the rational homotopy type. \ With these general explanations
out of the way, let us return to our twisted $K$-theory over $\R.$ \ Let $\eta$
 denote a fixed closed 3-form on $X$ representing the class of $P.$ \ Then
a class in $E_{4}$ in the spectral sequence may be represented by a closed
$p$-form $x$ with $\left[  \eta\right]  \left[  x\right]  =0.$ Since $\left[
\eta\right]  ^{2}=0$, we can form the Massey triple product $\left\{
\eta,\eta,x\right\}  .$ \ Since we actually have $\eta^{2}=0$ as a form, we can take
\begin{equation}
\left\{  \eta,\eta,x\right\}  =v\eta\;\;\text{where }\eta x=dv,\tag{5.2}%
\end{equation}
and the Massey product is less ambiguous than in the general case.
If we change the choice of $v$ by adding $v_{0}$ with $dv_{0}=0,$ then $\eta
v_{0}$ belongs to the boundaries $B_{4}$ and so its class in $E_{4}=E_{5}$ is
well-defined. \ Thus we have exhibited an operation
\[
\delta_{5}:E_{4}^{p}\longrightarrow E_{4}^{p+5}%
\]
In \S 7 we shall show that $-\delta_5$ is the differential $d_5$ of the spectral seqence for twisted $K$-theory.

\bigskip

 If $\{\eta ,\eta ,x \} =0$ we can repeat the process, defining a Massey
product
\[
\left\{  \eta,\eta,\eta,x\right\}  =\eta w\;\;\text{where\ \ }\left\{
\eta,\eta,x\right\}  \equiv dw\ (\operatorname{mod}\ \text{multiples of }\eta),
\]
giving an operation
\[
\delta_{7}=E_{6}^{p}\longrightarrow E_{6}^{p+7}%
\]
which we shall see in \S 7 is minus the differential $d_7$ of the spectral sequence.

\bigskip

 Continuing in this way, we find arbitrarily many potentially non-vanishing differentials. \ This
means that, in general, $K_{P}^{\ast}\left(  X\right)  $ will have much
smaller rank than $H^{\ast}\left(  X\right)  .$ \ A more precise formulation
will emerge in the subsequent sections. \ However, one fact can already be
noted. \ If, for a given manifold $X,$ all Massey products vanish, then the
spectral sequence stops at $E_{4}$, so that (5.1) holds and 
\begin{equation}
{\rm gr}\; K_{P}^{\ast}\left(  X\right)  \cong E_{4}\tag{5.3}%
\end{equation}
An important class of manifolds for which this holds are \emph{compact
K\"{a}hler manifolds}, where a deep theorem [DGMS] asserts that all
Massey products vanish. On the other hand we shall give examples in the Appendix of manifolds for which the Massey products of arbitrarily high order do not vanish. 

\bigskip

 In ordinary $K$-theory the collapsing of the Atiyah-Hirzebruch
spectral sequence over $\Q$ (or $\R)$ is intimately tied up with the existence
of the \emph{Chern character}, which gives a ring isomorphism
\[
{\rm ch}:K^{\ast}\left(  X\right)  \otimes \Q\longrightarrow H^{\ast}\left(
X,\Q\right)  .
\]

 As we have seen, the twisted K-theory given by a twisting element
$\eta$ of infinite order is far from being isomorphic over $\Q$ to cohomology:
typically it is considerably smaller. \ So the question arises: \ is there an
appropriate twisted cohomology theory with a corresponding twisted Chern
character? \ In the next section we shall describe such a theory.

\section{Twisted cohomology}

 \ \ \ \ Given a manifold $X$ and a closed 3-form $\eta$ we shall introduce a
\textquotedblleft twisted\textquotedblright\ cohomology theory (over $\R$).
\ On the de Rham complex $\Omega (X)$ of differential forms on $X$ we define an operator $D_{\eta}$ (or just $D$)
\begin{equation}
D=d-\eta\tag{6.1}%
\end{equation}
Note that
\[
D^{2}=\left(  d-\eta\right)  \left(  d-\eta\right)  =d^{2}-d\eta-\eta
d+\eta^{2}=0,
\]
since $\eta$ has odd degree. \ However, $D$ is not homogeneous, though it
preserves the grading $\operatorname{mod}2$. \ Thus the cohomology groups of
$D$ can be defined and they are graded $\operatorname{mod}2$. \ We denote
these groups by
\[
H_{\eta}^{0}\left(  X\right)  \;\;,\;\;H_{\eta}^{1}\left(  X\right)  .
\]

If $Y$ is a submanifold of $X$ we can of course define relative groups $H^*_{\eta}(X,Y)$ as the $D$-cohomology of the relative de Rham complex $\Omega (X,Y)$ of forms which vanish on $Y$; and we shall have the usual long exact sequence.
\bigskip

\no{\bf Remarks}

\bigskip

(i)\ \ It is not clear how one could define twisted cohomology with {\it integer} coefficients. Even on the level of cohomology we need not have $\eta^2 =0$ integrally (though, as we saw in Proposition 4.6, we do have $(Sq^3_{\Z}-\eta)^2 =0$).

\bigskip

(ii)\ \ The groups $H^*_{\eta}(X)$ depend on the closed form $\eta$ and not just on its cohomology class. If $\eta$ and $\eta '$ are cohomologous --- say $\eta ' - \eta = d\zeta$ --- then multiplication by the even-dimensional form $e^{\zeta}$ induces an isomorphism $H^*_{\eta}(X)\to H^*_{\eta '}(X)$, as $D_{\eta '}\circ e^{\zeta} = e^{\zeta}\circ D_{\eta}$. \ Because we have  $(d\alpha)\beta = D_{\eta}(\alpha \beta)$ if $D_{\eta}\beta = 0$ \ the isomorphism of twisted cohomology evidently does not change if $\zeta$ is altered by the addition of an exact 2-form.  \ One consequence of this --- taking $\zeta$ to be closed --- is that the group $H^2(X;\R)$ acts on the twisted cohomology $H^*_{\eta}(X)$. \ We shall see in the next section that this action corresponds by the Chern character to the natural action of the Picard group $H^2(X;\Z)$ on $K^*_P(X)$ which was described in [AS]. \ Another consequence is that if we are given not just one closed 3-form $\eta$ but a family $\{\eta_a\}_{a\in A}$ of 3-forms together with a choice of a 2-form $\zeta_{ab}$ for each $a,b\in A$ such that $\eta_b - \eta_a = d\zeta_{ab}$, and if
\[\zeta_{ab}-\zeta_{ac} +\zeta_{bc} = 0 \]
for each $a,b,c$, then the groups $H^*_{\eta_a}(X)$ are {\it canonically} isomorphic. \ We shall sometimes write them as $H^*_{\F}(X)$, where $\F$ denotes the family $\{\{\eta_a\},\{\zeta_{ab}\}\}.$ \ In fact   the forms $\zeta_{ab}$ need to be given only up to the addition of exact forms. \ We shall refer to such a family $\F$ as a {\it coherent family} of closed 3-forms.  

\bigskip

(iii)\ \ We shall often encounter the situation where we have an open covering $\mathcal U = \{U_a\}$ of $X$ and a closed 2-form $\omega_{ab}$ defined in each intersection $U_{ab}=U_a\cap U_b$ such that
\[\omega_{ab} -\omega_{ac} +\omega_{bc} = 0\]
in the triple intersection $U_{abc}$. \ 
By choosing a partition of unity $\{\lambda_a\}$ subordinate to the covering $\mathcal U$ we can then define a global closed 3-form $\eta$ whose value in $U_a$ is $dB_a$, where $B_a =\sum_b \lambda_b.\omega_{ab}$. (Indeed, in string theory it is usually the ``B-field" represented locally by the $B_a$ which is given, and $\omega_{ab}$ is defined as $B_a -B_b$.) \ The previous remark shows that the mod 2 graded de Rham complex $\Omega(X)$ of $X$ with the twisted differential $D_{\eta}$ is isomorphic to the complex $\Omega_{\omega}(X)$ in which a form is a collection $\alpha_a \in \Omega(U_a)$ such that $\alpha_b = e^{\omega_{ab}}\alpha_a$ for all $a,b$, with the usual exterior differential.  \ The map takes the collection $\alpha_a$ to the globally defined form which is $e^{B_a}\alpha_a$ in $U_a$. \ (Furthermore, if the covering $\mathcal U$ and the forms $\omega_{ab}$ are given, then varying the partition of unity gives rise to a coherent family $\F$ of 3-forms in the sense of the previous remark.)

\bigskip
 
\bigskip

\bigskip

 Although the twisted differential $D$ does not preserve the grading of the de Rham complex, it does preserve the
filtration whose $p$th stage is the sum of the forms of all degrees $\geq p$. \ The filtration therefore gives us a spectral sequence converging to the twisted cohomology. \ The
$E_{1}$ term is clearly just the differential forms with the usual $d$
(since $\eta$ raises filtration by $3$), and so
\[
E_{2}^{p}=H^{\ast}(X,\R)
\]
is the usual de Rham cohomology. \ To compute the higher differentials we
proceed as follows. \ Let
\[
x=x_{p}+x_{p+2}+x_{p+4}+...
\]
be an (inhomogenous) form with $dx_{p}=0,$ so that $x_{p}$ represents a class
$\left[  x_{p}\right]  \in H^{p}.$ \ Then
\[
Dx=\left(  d-\eta\right)  x=dx_{p}+\left(  dx_{p+2}-\eta x_{p}\right)
+\left(  dx_{p+4}-\eta x_{p+2}\right)  +...
\]
i.e
\[
y=Dx=y_{p+1}+y_{p+3}+y_{p+5}...
\]
where
\begin{align*}
y_{p+1} &  =dx_{p}=0\\
y_{p+3} &  =dx_{p+2}-\eta x_{p}\\
y_{p+5} &  =dx_{p+4}-\eta x_{p+2}\\
&  ....................
\end{align*}
The class $\left[  y_{p+3}\right]  =-\left[  \eta\right]  \left[  x_{p}\right]
$ represents $d_{3}x,$ so that the differential $d_{3}$ of the spectral
sequence is just multiplication by $-\eta.$ \ Hence $E_{4}$ is the same as the
$E_{4}$ term of the Atiyah-Hirzebruch spectral sequence (over $\R$) that we
found in Section 4.

\bigskip

 Proceeding further, if $\left[  \eta x_{p}\right]  =0$ so that $\eta
x_{p}=dv,$ we can choose $x_{p+2}=v$ and make $y_{p+3}=0.$ \ Then
$y_{p+5}=dx_{p+4}-\eta v$ represents $d_{5}x_{p}.$ \ In other words
\[
d_{5}\left[  x_{p}\right]  =-\left\{  \eta,\eta,x_{p}\right\}
\]
is given by the Massey triple product. \ Continuing this way we see that
\[
d_{7}\left[  x_{p}\right]  =-\left\{  \eta,\eta,\eta,u\right\}
\]
and so on. \ Thus we have proved
\bigskip

\noindent{\bf Proposition 6.1} \ \ {\it The iterated Massey products with }$\eta$
 {\it give (up to sign) all the higher differentials of the spectral sequence for the
twisted cohomology.}

\bigskip

 The similarity of this spectral sequence with that of twisted
$K$-theory described in Section 5 is clear. \ In the next section we shall bring
them together via the \emph{twisted Chern character} which we shall define.

\bigskip

Let us first, however, point out a simple corollary of the existence of the filtration spectral sequence which makes working with twisted cohomology a lot more flexible. Obviously we can  define twisted cohomology $H^*_{\eta}(A)$ whenever we have a differential graded algebra $A$ which is graded-commutative, and a chosen closed element $\eta$ in degree 3. The spectral sequence tells us that if $\phi : A \to A'$ is a homomorphism of differential graded algebras which induces an isomorphism of cohomology in the usual sense, then it induces an isomorphism 
\[H^*_{\eta}(A) \to H^*_{\phi (\eta )}(A') \]
for every closed $\eta \in A^3.$  \  Furthermore, if instead of $\eta$ we have a closed 3-form $\eta '$ in $A'$, then we can define a coherent family $\phi ^*(\eta ')$ in $A$ such that
\[ H^*_{\phi ^*(\eta ')}(A) \cong H^*_{\eta '}(A'). \]
The family is indexed by the set of pairs \  $(\eta , \xi ')$ \ in $A^3 \oplus (A')^2$ such that $d\eta = 0$ and $\phi (\eta ) +d\xi ' = \eta '$:\  this gives us a coherent family because the mapping-cone\footnote{The mapping-cone of $A \to A'$ is the total complex of the double complex $\ldots \to 0 \to A \to A' \to 0 \to \ldots \ .$} of the equivalence $A \to A'$ is acyclic.

An important application of this remark is the following. \ Any topological space $X$ has a ``singular de Rham complex" \  $\Omega_{\mathrm{sing}}(X)$, \ in which a $p$-form is a family of $p$-forms \ $\alpha _{\sigma} \in \Omega^p(\Delta^m)$, \  one for each singular simplex \  $\sigma :\Delta ^m \to X$, \  which are compatible under the face and degeneracy maps. \ It is well-known [Su] that the cohomology of \ $\Omega_{\mathrm{sing}}(X)$ --- which is clearly a commutative differential graded algebra --- is the singular cohomology $H^*(X;\R)$. \ Now if $X$ is a smooth manifold there are natural homomorphisms of differential graded algebras
\[\Omega (X) \to \Omega_{\mathrm{sm-sing}}(X) \leftarrow \Omega_{\mathrm{sing}}(X) \]
inducing isomorphisms of cohomology, where \ $\Omega_{\mathrm{sm-sing}}(X)$ \  is defined like \ $\Omega_{\mathrm{sing}}(X)$ \  but using only {\it smooth} singular simplexes. \ Coherent families of 3-forms in these algebras correspond one-to-one, and the same twisted cohomology is obtained from any of them. \ In other words, we need not hesitate to speak of the twisted cohomology of spaces which are not manifolds.

\bigskip

We end this section with a technical result needed in the next section. \ Let us first notice that if $X$ is a space for which $H^{{\rm odd}}(X;\R) = 0$ then the spectral sequence shows that $H^{{\rm odd}}_{\eta}(X) =0$ for any $\eta$ \ ; \  and the same principle holds for relative groups. \ We shall apply this to the situation studied in section 3, where we have a fibration $Y\to S^3$ over the 3-sphere, with fibres $F_0$ and $F_{\infty}$ over the poles 0 and $\infty$ in $S^3$. \ We are interested in the twisted cohomology $H_{\eta}(Y)$, where $\eta \in \Omega^3(Y)$ is pulled back from a form on $S^3$ which generates $H^3(S^3)$, but vanishes near the base-point 0. \ Thus $\eta$ vanishes in the neighbourhood of $F_0$, and so there is a restriction map $H^{{\rm ev}}_{\eta}(Y)\to H^{{\rm ev}}(F_0)$.

\bigskip

\noindent{\bf Proposition 6.2} \ {\it For the fibration} $Y \to S^3$ {\it with fibre} $F_0$, {\it if} $H^{{\rm odd}}(F_0 ;\R) =0$ {\it then}
$$H^{{\rm ev}}_{\eta}(Y)\to H^{{\rm ev}}(F_0)$$
{\it is injective.}

\bigskip

\noindent{\bf Proof} \ \ As explained in \S 3, we have $H^i(Y,F_0)\cong H^{i-3}(F_{\infty})$, so the hypothesis implies that the relative group $H^{{\rm ev}}_{\eta}(Y,F_0)$ vanishes. \ \ $\Box$

\section{The Chern character}

\bigskip 

 \ \ \ \ For each projective bundle $P$ on $X$ we should like to have a natural Chern character map

\[ {\rm ch}:K^*_P(X) \to H^*_{\eta}(X)\]
which induces an isomorphism $K^*_P(X)\otimes \R \to H^*_{\eta}(X)$. \ Here $\eta$ is a 3-form on $X$ representing the class of the projective bundle $P$. \ In fact it is unreasonable to try to assign a specific 3-form to the bundle $P$ without choosing any additional structure\footnote{What is needed is a choice of a {\it string connection} in the sense of [S3], or of some version of a {\it gerbe connection} [BCMMS].}, and we shall instead define a coherent family of forms in the sense explained in Remark (ii) of the previous section.

\bigskip

In this paper we are using the traditional methods of algebraic topology --- though we have compromised a little by defining twisted cohomology in terms of the de Rham complex --- and in that spirit we shall define the Chern character as a universal form on the representing space for twisted $K$-theory.  \ Before doing so, however, it seems worthwhile to sketch briefly how things would look from the point of view of Chern-Weil theory, for that is more geometrical, and more relevant in applications of the theory. \ To treat the Chern-Weil theory fully would require a careful discussion of smooth Hilbert bundles and families of unbounded operators in them, and we shall not embark on that.

\bigskip

The Chern character of a smooth vector bundle $E$ on $X$ is represented by the inhomogeneous differential form tr$(e^{F_E})$, where $F_E \in \Omega^2(X;{\rm End}(E))$ is the curvature of an arbitrary connection in $E$. \ (We recall that if 
\[d_E :\Omega^i(X;E)\to \Omega^{i+1}(X;E)\] is the covariant exterior derivative operator associated to a connection in $E$ then $d_E^2$ is multiplication by $F_E$: a little cryptically, we write $d_E^2 = F_E$.)
\bigskip

It is perhaps slightly less well known that there is a generalization --- due to Quillen [Q] ---  of this description of the Chern character to the case of an element of $K(X)$ represented by a family of unbounded self-adjoint degree 1 Fredholm operators $\{T_x\}$ in a mod 2 graded bundle $\H$ of Hilbert spaces on $X$, at least when $e^{-T_x^2}$ is of trace class for all $x\in X$. \ We simply form the {\it superconnection} $d_{\H} + {\rm i}T$, an operator on $\Omega(X;\H)$, and define its curvature as 
\[ \F = (d_{\H} + {\rm i}T)^2. \]
Then the inhomogeneous differential form\footnote{The {\it supertrace} of an operator $T$ in a mod 2 graded vector space $\H = \H^0 \oplus \H^1$ is defined by
\[ {\rm str}\;T={\rm tr}\;T|\H^0 -{\rm tr}\;T|\H^1. \]} str$(e^{\F})$ represents the Chern character. \ (We are assuming here that the bundle $\H$ is smooth in the sense that one can speak of differential forms on $X$ with values in $\H$, that we have a connection in $\H$ giving an operation $d_{\H}$ of covariant differentiation on such forms, and that the family $\{T_x\}$ is appropriately smooth. \ For more details, see [BGV].)
\bigskip

Quillen's description of the Chern character is well adapted to the definition of twisted $K$-theory  $K_P(X)$ by means of families of Fredholm operators in a projective Hilbert bundle $P$. \ It is simplest to explain this in terms of local trivializations. Let us choose local trivializations of $P$ over the sets of an open covering $\mathcal U = \{U_a\}$ as above, and represent $P$ and its connection by an actual Hilbert bundle $\H_a$ with a unitary connection $d_{\H_a}$ in $U_a$.  \ Then over $U_{ab}=U_a \cap U_b$ we have
\[\H_b \cong \H_a \otimes L_{ab},\]
where $L_{ab}$ is a line bundle with a unitary connection. \ Let the curvature of $L_{ab}$ be the closed 2-form $\omega_{ab}$ on $U_{ab}$.  \ Then the curvature forms $\F_a$ defined in each $U_a$ are related by
\[\F_b = \F_a + \omega_{ab}\]
in $U_{ab}$, and so the character forms str$(e^{\F_a})$ constitute an element of the twisted de Rham complex $\Omega_{\omega}(X)$ which --- as we have pointed out in Remark (iii) in \S 6 --- calculates the twisted cohomology $H^*_{\eta}(X)$. \ If the twisting is given by a B-field $\{B_a\}$ then the Chern character is given by the globally defined form str$(e^{\F_a + B_a})$. \ This is the Chern-Weil definition of the twisted Chern character.

\bigskip

Now let us give a less problematic definition of the Chern character as a twisted cohomology class of the universal fibration Fred$(P)$ of Section 3, with the twisting defined by the generator of $H^3(K(\Z,3);\Z)$. \ The first point to notice is that, by the remark at the end on Section 6, we can replace the universal fibration by its restriction $Y$ to the sphere $S^3$ contained in the base-space $K(\Z,3)$, for the two spaces are rationally homotopy equivalent, and so the restriction map of their singular de Rham complexes induces an isomorphism of both ordinary and twisted cohomology. \ As in Section 3 we shall denote the fibres of the universal bundle over the north and south poles of $S^3$ by $F_{\infty}$ and $F_0$. \  We shall prove that there is a {\it unique} twisted class of $Y$ which restricts to the standard Chern character on the fibre $F_0$. \ (We shall explain the significance of the uniqueness presently.) \ Since we are concentrating on the case when the twisting class
$\eta$ is of infinite order we may (as explained in Section 3) restrict
ourselves to the index zero component of Fred$(P)$. \ Let $s_{n}$ be a closed form on  $F_0$ representing the cohomology generator  denoted by
the same symbol in Section 3, and let $\eta$ be a 3-form on $Y$ pulled back from a generating 3-form on $S^3$ whose support does not contain either pole of $S^3$. \ Our task is to construct a sequence of closed forms $S_n$ on $Y$ with the two properties
\bigskip

(i)\; $S_n |F_0 = s_n $ ,
\bigskip

\noindent and
\bigskip

(ii)\;  $dS_n =  n \eta S_{n-1}$.
\bigskip

We shall do this by induction on $n$. \ Suppose that $S_m$ has been found for $m<n$. \ Because the bundle $Y$ can be trivialized in the complement of the single fibre $F_{\infty}$ we can pull back $s_n$ to a closed form on this complement, and then extend it further to a not-necessarily-closed  form $S'_n$ on all of   $Y$. \ (This extension is especially easy using the singular de Rham complex, as we can extend by zero; if we were working with honest smooth forms we should have to multiply the pull-back by a suitable real-valued function equal to 1 outside a neighbourhood of $F_{\infty}$ and vanishing in a smaller neighbourhood of $F_{\infty}$.) \  Then $dS'_{n}=0$ away from $F_{\infty}$, and so it defines
a class in $H^{2n+1}\left(  Y,F_{0}\right)  .$ \ By the arguments in Section 3
(leading up to Proposition (3.9)) this class can also be represented by  the form $ n\eta S_{n-1}$, which also vanishes near $F_0$.
 \ Hence on the level of forms we have
\begin{equation}
dS'_{n}-n\eta S_{n-1}=d\theta_{n},\tag{7.1}%
\end{equation}
where $\theta_{n}$ is form on $Y$ with support away from $F_{0}.$ \ So we define $S_n = S'_n - \theta_n$.

\bigskip

\ Note that the formulae (i) and (ii) still hold for $n=1$, where $S_{0}=0$ in the zero-component.

\bigskip

\noindent Thus the even form on $Y$%
\begin{equation}
{\rm ch}=\sum\limits_{n=1}^{\infty}\dfrac{S_{n}}{n!}\tag{7.2}%
\end{equation}
satisfies $\left(  d-\eta\right)  {\rm ch}=0,$ showing that it is a cocyle for the
twisted cohomology of the universal fibration.

\bigskip

 The construction of the form in (7.2)  provides a class in $H_{\eta}\left(  Y\right)  $ which lifts the usual
Chern character in $H^*\left(  F_{0}\right)  .$   \ It defines a twisted Chern character in the sense that when we have a manifold $X$ with a projective bundle $P$ we can choose a smooth map  $f:X\to K(\Z,3)$ which pulls back $P$ from the universal $\C P^{\infty}$-bundle $\hat{P}$ on $K(\Z,3)$. \  Then $K_P(X)$ is defined as the set of homotopy classes of lifts $F$ of $f$ to $Y$, and for each such lift $F$ we define  ${\rm ch}\,F =\,F^*({\rm ch})$.\  This is an element of the twisted cohomology $H^*_{\eta _f}(X)$, where $\eta_f = f^*(\eta)$. \  The definition of the twisting form $\eta_f$ depends on the choice of the classifying map $f$, but as $f$ varies we obtain a coherent family $\F_P$ of forms in the sense of the preceding section. \  For though the space of maps $f$ pulling back a bundle equivalent to $P$ is not contractible, the space of pull-back diagrams
\[\begin{array}{rcc}P&\longrightarrow &  \hat{P}\\
\downarrow &&\downarrow\\
X&\longrightarrow & K(\Z,3),
\end{array} \]\
i.e. of {\it pairs} $(f,\alpha)$, where $\alpha$ is an isomorphism $\alpha :P\to f^*\hat{P}$,\, {\it is} contractible, and we index $\F_P$ by this  set of pairs. \ Then  there is a path, unique up to homotopy, from one choice of $f$ to any other --- say $f'$ ---, and hence, by the usual argument of de Rham cohomology, a choice of $\zeta_{ff'}$, unique up to an exact form, such that 
\[d\zeta_{ff'}= \eta_f -\eta_{f'} .\]

\bigskip

The uniqueness of the twisted class of $Y$ which lifts the Chern character on $F_0$ is important.\footnote{We are indebted to U.Bunke for a helpful discussion of this point.} \ The classical Chern character defines a ring homomorphism from $K$-theory to cohomology, but if we wanted just an additive homomorphism we could equally well use $\sum \lambda _k \ {\rm ch}_k$ for any sequence of rational numbers $\{\lambda_k\}$. \ The position is different with twisted $K$-theory. \ A class $\xi \in H^*_{\eta}(Y)$ defines an additive transformation from twisted $K$ to twisted cohomology if and only if it is {\it primitive} (i.e. maps to $\xi \otimes 1 + 1 \otimes \xi)$ under the fibrewise Whitney sum map
\[ Y\times_B Y \to Y, \]
where $B=K(\Z,3)$. \ Now Proposition 6.2 applies to $Y\times _B Y$, and it  follows that a twisted class $\xi$ defines an additive transformation if and only if its restriction to the fibre is primitive, and hence of the form $\sum \lambda_k \ {\rm ch}_k$. \ But only if all the $\lambda_k$ are equal does this class lift to a twisted class. \ Thus the twisted Chern character is singled out by additivity alone, without mentioning any multiplicativity property. \ On the other hand, the argument just given also proves that our twisted character is indeed multiplicative in the natural sense with respect to the pairings $$K_P(X)\otimes K_{P'}(X) \to K_{P\otimes P'}(X)$$ and $$H_{\eta}(X)\otimes H_{\eta '}(X) \to H_{\eta + \eta '}(X).$$
\bigskip   

We have now defined the twisted character on the even group $K_P(X)$. \  But we can extend the definition to a homomorphism
\begin{equation}
{\rm ch}:K_{P}^{\ast}\left(  X\right)  \longrightarrow H_{\F_P}^{\ast}(X)\nonumber
\end{equation}
of $\operatorname{mod}2$ graded theories simply by replacing $X$ by $X\times S^1$ in the usual way.

\bigskip

The remaining formal properties of the twisted Chern character also follow from
those of the untwisted case, as we shall now explain.

\bigskip

 Because the twisted Chern character is functorial it is compatible with
the filtrations which are used to construct the Atiyah-Hirzebruch spectral sequences, and so it induces a homomorphism of the associated spectral
sequences of Sections 4 and 6. \ But when the $K$-theory spectral sequence is tensored with $\R$ the $E_{2}$ terms coincide, and
the homomorphism there is an isomorphism. \ Hence the whole spectral sequence
map is an isomorphism, and we conclude that
\[
{\rm ch}:K_{P}^{\ast}\left(  X\right)  \otimes \R\to H_{\F_P}^{\ast}\left(
X\right)
\]
is an isomorphism. \ This argument also shows that the differentials
of the spectral sequences must agree. \ Since the higher differentials for
twisted cohomology were shown in Section 6 to be given precisely by the
iterated Massey products, the same follows for the twisted $K$-theory spectral
sequence. 

\bigskip

 To summarize, we have established\bigskip

\noindent\textbf{Proposition (7.4) \ }\emph{There is a functorial twisted Chern character
}
\[
{\rm ch}:K_{P}^{\ast}\left(  X\right)  \longrightarrow H_{\F_P}^{\ast}\left(
X\right)
\]
\emph{which becomes an isomorphism after tensoring with }$\R.\bigskip$

\noindent\textbf{Proposition (7.5) \ }\emph{In the twisted Atiyah-Hirzebruch
spectral sequence over }$\R$\emph{\ the higher differentials }$d_{5,}d_{7}%
,...$\emph{\ are given by the iterated Massey products }
\begin{align*}
d_{5}\left(  x\right)   & =-\left\{  \eta,\eta,x\right\} \\
d_{7}\left(  x\right)   & =-\left\{  \eta,\eta,\eta,x\right\} \\
& ...
\end{align*}

\section{Chern Classes}

 \ \ \ \ Ordinary (untwisted) $K$-theory is related to cohomology in three
different ways: via the Atiyah-Hirzebruch spectral sequence, via the Chern
character, and via Chern classes. \ These are all connected in well-known ways.
\ In particular the Chern character can be expressed in terms of the Chern
classes, and the converse holds over the rationals. \ However, in the twisted
case (for a twisting element of infinite order) there are more substantial
differences, and one cannot pass so easily from one to another. \ For example,
the twisted Chern character, not being graded, does not lead to a sequence of rational
twisted Chern classes.\bigskip

 The basic result about characteristic classes in the untwisted case is that any cohomology class naturally associated to an element $\xi \in K(X)$ is a polynomial in its Chern classes $c_n(\xi)$. \ For a twisted class $\xi \in K_P(X)$ we can look for characteristic classes in either the ordinary cohomology or the twisted cohomology of $X$.  \ Although the latter question might seem
more logical, we shall begin with the former, which is easier and more
tractable. \ Our results are summarized in Proposition 8.8 below.\bigskip

 For any generalized cohomology theory characteristic classes correspond to cohomology classes of the classifying space. \ For twisted $K$-theory
the appropriate classifying space is the universal fibration studied in
Section 3, with fibre Fred$(\H)$ and base the Eilenberg-MacLane space $K\left(
\Z,3\right)  .$\bigskip

 The cohomology classes of the total space define characteristic
classes for our twisted $K$-theory. \ For a twisting by an element $\eta$ of
infinite order only the index zero component of Fred$\left(  \H\right)  $ is
relevant, and so we shall restrict ourselves to this case.\bigskip

 We now consider only real or rational coefficients, so that, as in
Section 3, we can replace the Eilenberg-MacLane space by the 3-sphere $S^{3}$.
\ The key result of Section 3 was Proposition (3.9), determining the only
non-trivial differential $d_{3}$ for the Leray spectral sequence of the
fibration. \ As noted there, the Leray spectral sequence 
reduces to the Wang exact sequence (3.5) , and hence gives the
exact sequence:
\begin{equation}
0\rightarrow H^{\text{ev}}\left(  Y\right)  \overset{j}{\rightarrow
}H^{\text{ev}}\left(  F\right)  \overset{d}{\rightarrow}H^{\text{ev}}\left(
F\right)  \rightarrow H^{\text{odd}}\left(  Y\right)  \rightarrow0,\tag{8.1}%
\end{equation}
where we have lumped together the even and odd degrees. \ Here $F=\ $%
Fred$_{0}\left(  \H\right)  $ and $Y$ is the total space.\bigskip

 The homomorphism $d$ is the derivation given by (3.9) with $k=0$:
\begin{align}
d\;s_{n}  & =n\;s_{n-1}\;\;\;\;\;\;\;\;n\geq2\nonumber\\
d\;s_{1}  & =0\tag{8.2}%
\end{align}
and (8.1) shows that
\begin{equation}
H^{\text{ev}}\left(  Y\right)  =J,\tag{8.3}%
\end{equation}
where $J\subset Q\left[  s_{1},s_{2},...\right]  $ is the kernel of $d$ in
(8.2).\bigskip\ \  The sequence (8.1) also determines $H^{\text{odd}}\left(  Y\right)  $ as
the cokernel of $d$, but, as the following lemma shows, this is zero except
for $H^{3}.$
\bigskip

\noindent{\bf Lemma 8.4} \ \ {\it The derivation $d$ given by (8.2) is surjective in degree $>0.$}
\bigskip

\noindent{\bf Proof}\ \ For simplicity of notation it is convenient to
introduce new variables $x_{n}\left(  n\geq 0\right)  $ by
\[
x_{n}= {\rm ch}_n =\dfrac{s_{n}}{n!}\;\;\;\;\;\;\;\;\;n\geq1,
\]
and to put $x_{0}=0.$ Then $d$ acting on a polynomial $f\left(  x_{1}%
,...\right)  $ becomes
\[
df=\sum\limits_{n\geq1}\dfrac{\partial f}{\partial x_{n}}x_{n-1}%
\]
\bigskip

\noindent It is sufficient to consider monomials of degree $n>0$, and we
assume, as an inductive hypothesis, that the Lemma is true for all smaller
positive values of $n.$ \ If $f$ has degree $n$ and is not in the image of
$d,$ then we can write $f=x_{k}g$ with $k$ minimal but $k<n$ (since
$x_{k}=dx_{k+1}).$ \ Then $0<\deg\left(  g\right)  <n$, so by the inductive
hypothesis $g=dh$ for some $h.$ \ Then
\[
f=x_{k}g=x_{k}dh=d\left(  x_{k}h\right)  -x_{k-1}h
\]
so that $x_{k-1}h$ is also not in the image of $d.$ \ But, by the minimality
of $k,$ this contradicts our assumption, proving the lemma. $\Box$
\bigskip

 Thus $H^{q}\left(  Y\right)  =0$ for odd $q\neq 3,$ while 
$H^{3}\left(  Y\right)  \cong H^{3}\left(  S^{3}\right)  $ with the generator
pulled back from the base. \ Thus, apart from the twisting class $\eta$ itself, \emph{there are no odd
(rational) characteristic classes for twisted }$K$\emph{-theory.}\bigskip

 The even characteristic classes are then given by the subring $J.$
\ We shall describe some aspects of this ring, but it does not seem easy to
give its full structure in terms of generators and relations.\bigskip

 One consequence of (8.1) and (8.4) is a formula for the Poincar\'{e}
series of $J.$ \ Let $A=\underset{n}{\oplus}A_{n}=\Q\left[  x_{1}%
,x_{2},...\right]  $ be the graded polynomial algebra, where degree
$x_{n}=n,$ and let
\[
J=\oplus J_{n}%
\]
be the grading of $J$. \ (Note that $d$ in (8.1) just shifts the grading by 1.)
\ Then (8.1) and (8.4) lead to the short exact sequence
\[
0\rightarrow J_{n}\rightarrow A_{n}\overset{d}{\rightarrow}A_{n-1}%
\rightarrow0
\]
for $n\geq2$, while
\[
J_{0}=A_{0}\;\;,\;\;J_{1}=A_{1}.
\]
Put $a_{n}=\dim A_{n},\;\;j_{n}=\dim J_{n}.$ \ Then
\begin{align*}
j_{n}  & =a_{n}-a_{n-1}\;\;\;\;\;\left(  n\geq2\right) \\
j_{1}  & =a_{1}=1\\
j_{0}  & =a_{0}=1.
\end{align*}
Defining the Poincar\'{e} series $J\left(  t\right)  ,\;A\left(  t\right)  $
by
\[
J\left(  t\right)  =\sum\limits_{n\geq0}j_{n}t^{n}\;\;\;\;\;\;A\left(
t\right)  =\sum\limits_{n\geq0}a_{n}t^{n},%
\]
we deduce that
\[
\sum\limits_{2}^{\infty}j_{n}t^{n}=\sum\limits_{2}^{\infty}a_{n}t^{n}%
-\sum\limits_{2}^{\infty}a_{n-1}t^{n},%
\]
or
\begin{align*}
J\left(  t\right)  -\left(  1+t\right)   & =A\left(  t\right)  -\left(
1+t\right)  -t\left[  A\left(  t\right)  -1\right] \\
J\left(  t\right)   & =A\left(  t\right)  -tA\left(  t\right)  +t\\
& =\left(  1-t\right)  A\left(  t\right)  +t.
\end{align*}
But, since $A$ is the polynomial algebra,
\[
A\left(  t\right)  =\dfrac{1}{\left(  1-t\right)  \left(  1-t^{2}\right)  ...},%
\]
and we have proved
\begin{prop}
The Poincar\'{e} series of $J$ is given by
\begin{equation}
J\left(  t\right)  =\dfrac{1}{\left(  1-t^{2}\right)  \left(  1-t^{3}\right)
...}+t.\tag{8.5}%
\end{equation}
\end{prop}

 The first few explicit values of $j_{n}$ are
\begin{align*}
j_{0}  & =1,\;j_{1}=1,\;j_{2}=1,\;j_{3}=1,\;j_{4}=2,\;j_{5}=2,\\
j_{6}  & =4,\;j_{7}=4,\;j_{8}=7.
\end{align*}

This shows that $J$ is not just a polynomial algebra (on generators of degrees
$\geq2$), as one might have expected from comparison with the
finite-dimensional groups $U\left(  N\right)  \rightarrow PU(N)$. \ The extra
term $t$ in (8.5) arises from the class $x_{1}=s_{1},$ but $J$ does not even contain a polynomial subalgebra on generators of degrees $\geq2$, for $x_1^2$ and $x_1^3$ are the only elements in degrees 2 and 3, and they are algebraically dependent.
\bigskip

 Although $J\subset A$ is explicitly defined by the formulae (8.2)
there is a more fundamental way of describing it that was already implicit in
Section 2, where we studied the action of tensoring by line-bundles. \ We
continue with our present notation, where we have replaced the power sums
$s_{n}$ by the components $x_{n}$ of the Chern character. \ In terms of a
variable $u$ of degree $2,$ introduce variables $y_{1},y_{2},...$ related to
$x_{1},x_{2},...$ by
\begin{equation}
\sum y_{n}=\left(  \sum x_{n}\right)  \left(  \sum\dfrac{u^{k}}{k!}\right)
\tag{8.6}%
\end{equation}
Thus, explicitly
\begin{align*}
y_{1}  & =x_{1}\\
y_{2}  & =x_{2}+ux_{1}\\
y_{3}  & =x_{3}+ux_{2}+\dfrac{u^{2}}{2}x_{1}\\
& \ldots
\end{align*}
These equations are, of course, just a rewriting of equations 2.4 and 2.5.
\bigskip

Given a polynomial $f\left(  x_{1},x_{2}...\right)  $ in $A$, we define it to
be \emph{invariant }if
\begin{equation}
f\left(  y\right)  \equiv f\left(  x\right) \tag{8.7}%
\end{equation}
identically in $u.$ \ Now equation (8.6) is invertible (giving $x$ in terms of
$y)$, and it can be written symbolically as
\[
Y=X\;e^{u}%
\]
is invertible (giving $x_k$ in terms of
$y_1, \ldots ,y_k$). \  It defines an algebraic action of the additive group $\C$ (on which $u$ is the coordinate function) on $A$ by algebra automorphisms. \ This explains why we used the word \textquotedblleft
invariant\textquotedblright. \ (To keep track of the grading: if we think of the grading of $A$ as corresponding to an action of $\C^{\times}$ on it, then the action of $\C$ on $A$ is $\C^{\times}$-equivariant when $\C^{\times}$ acts by multiplication on $\C$.)
\bigskip

 As noted in Section 3, our homomorphism $d:A\rightarrow A$ is 
given by just using the linear term in $u$ in (8.6). \ In other words it is the
action of the Lie algebra of the 1-parameter group. \ But Lie algebra invariance implies Lie
group invariance. \ Hence the apparently much stronger identity (8.7),\, given
by using all higher powers of $u,$\,  is actually equivalent to the simpler one
using the linear terms. \ Thus \emph{our subring }$J\subset A$ \emph{is just
the invariant subring.}\bigskip

 A topological explanation of this fact is implicit in the formula (3.2) of Section 3. \ A class in $H^{\text{ev}}\left(  F\right)  $ which is in $J$
extends by definition to$\;H^{\text{ev}}\left(  Y\right)  ,$ where $Y$ is
fibred over $S^{3}.$ \ Because of the rational homotopy equivalence (3.3)  it
then also extends over $BPU\left(  H\right)  $, and in particular over the
subspace 
$\Sigma  \C P^{\infty}  $ of (3.2). \ But this passage from $\C P^{1}$ to
$ \C P^{\infty}$ is just the extension from the term linear in $u$ to all powers
of $u,$ recovering the algebraic argument above.\bigskip

 We summarize all this in
\bigskip

\noindent{\bf Proposition (8.8)} \ \ {\it Each polynomial in the universal Chern
classes which is invariant under the tensor product operation  (8.6) extends uniquely to a rational  characteristic class of twisted $K$-theory, and (apart from the twisting class) all such characteristic classes are obtained in this way.}

{\it Furthermore, all of the characteristic classes are annihilated by multiplication by the twisting class $\eta$.}
\bigskip

 Having identified the ring $J$ of characteristic classes of twisted
$K$-theory in \emph{ordinary} cohomology, we can now ask for characteristic
classes in \emph{twisted} cohomology. \ As explained in Section 6 we get the
twisted (real) cohomology as a $\operatorname{mod}2$ graded theory by starting
with the de Rham complex and using the total differential $D=d-\eta$ of (6.1).
\ We then get a spectral sequence, whose $E_{2}$-term is the usual real
cohomology and whose first differential $d_{3}$ is multiplication by the class
of $- \eta.$ \ The $E_{\infty}$ is then the associated graded group of our
twisted cohomology.\emph{\bigskip}

 We apply this to our universal space $Y$ . \ In this section we have shown
that
\[
H^{\text{ev}}\left(  Y\right)  =J\;\;\;\;\;\;H^{\text{odd}}\left(  Y\right)
=H^{3}\left(  Y\right)  =\R,
\]
and the only non-trivial differential $d_{3}$ gives an isomorphism
\[
H^{0}\left(  Y\right)  \rightarrow H^{3}\left(  Y\right)  .
\]
Thus the $E_{4}$-term of our spectral sequence is just the part
$J^{+}$ of $J$ of positive degree:
\begin{equation}
E_{4}\cong J^{+}.\tag{8.9}%
\end{equation}
Since there are no more terms of odd degree, all remaining differentials are
zero, and so
\[
J^{+}\cong E_{\infty}\cong {\rm gr}\;H_{\eta}\left(  Y\right),
\]
where $H_{\eta}$ is the universal twisted cohomology.\emph{\bigskip}

 To get characteristic twisted classes for twisted $K$-theory we
therefore have to choose liftings of $J^{+}$ back into the filtered ring
$H_{\eta}\left(  Y\right)  .$ \ For example, for the first class $x_{1}%
=s_{1}\in J^{2}$  the \emph{Chern character}, as constructed in Section
7, provides a natural lift. \ In general, a lift can be constructed by the following simple but not very illuminating formula.
{\bigskip}

Define a new derivation $\delta$ of the polynomial algebra $A=\C [x_1,x_2,\ldots]$ by 
\[\delta(x_k)=kx_{k+1} \]
for $k\geq 1$. \ The commutator $\Delta =[d,\delta]$ is then also a derivation: it multiplies each monomial in the $x_k$ by its degree, when each $x_k$ is given degree 1 rather than $k$. \ (This new grading of $A$ is preserved by both $d$ and $\delta$, and is hence inherited by $J$.) \ As $d$ and $\delta$ commute with $\Delta$ we have $[d,\exp (\lambda \delta)] = \lambda \Delta \exp (\lambda \delta)$ for any scalar $\lambda$, i.e.
\[(d-\lambda \Delta )\circ\exp (\lambda \delta ) = \exp (\lambda \delta )\circ d. \]
If $f$ is an invariant element of $A$ with $\Delta (f) = mf$ for some $m>0$, we accordingly have 
\[ (d-1)\exp (\lambda \delta )(f) = 0 \]
when $\lambda = 1/m$. \ This shows that $\exp (\lambda \delta )(f)$ defines a class in the twisted cohomology of the universal fibration which lifts $f$. When $f = s_1 = x_1$ it is precisely the Chern character.
\bigskip

 To find more natural lifts for  classes in $J^{+}$ one procedure
would be to first construct \emph{internal operations} in twisted $K$-theory
and then take their Chern characters. \ This is one motivation for our study
of such internal operators in Section 10. \ But, before we proceed down that
route, we will in the next section say something about an interesting subset
of elements of $J$ defined geometrically - the Koschorke
classes.\emph{\bigskip}

\section{Koschorke classes}

  \ \ \ \ The classical treatment of Chern classes has several aspects, all of
importance. \ One can introduce them as cohomology classes or perhaps
differential forms of the classifying space. \ Dually one can define them more
geometrically by explicit cycles on the classifying space - the famous
\textquotedblleft Schubert cycles\textquotedblright\  on the Grassmannian manifolds of subspaces of a vector space. \ These Schubert cycles
are defined relative to some fixed (partial) flag in the ambient vector space.
\ If we work with Fredholm operators in a Hilbert space we would need flags
formed by subspaces of finite codimension. \ But for a twisting of infinite order
 such flags do not exist for reasons explained in [AS]. \ However,
even though we cannot use the direct analogue of the Schubert cycles, there is
an alternative and very natural class of geometric cycles that we can use. \ These are the cycles studied by Koschorke [K].\emph{\bigskip}

 In Fred$\left(  \H\right)  $ we have the locally closed submanifolds
$F_{p,q}$ consisting of Fredholm operators $T$ with
\[
\dim\text{Ker\ }T=p\;\;\;\dim\;\operatorname{Coker}T=q.
\]
Of course $k=$ index $T=p-q$ is constant on the components Fred$_{k}\left(
\H\right)  ,$ so that for each $k$ we have a simple sequence of these
submanifolds.\emph{\bigskip}

 The submanifold $F_{p,q}$, or rather its closure $\bar{F}_{p,q}$,
defines an integral cohomology class $k_{p,q}$ of dimension $2pq$ called the
\emph{Koschorke class }.\emph{\bigskip}

 Clearly the submanifolds $F_{p,q}$ are invariant under the conjugation action of
$PU\left(  \H\right)  $ on Fred$\left(  \H\right)  $, and they therefore define a
submanifold of the total space of the associated fibration over $BPU\left(
\H\right)  ,$ and so over any pull-back to a space $X$ by a twisting class
$\eta.$ \ This exhibits a natural lift of the Koschorke classes in the fibre
Fred$\left(  \H\right)  $ to the total space, and hence to characteristic classes in the twisted cohomology $H_{\eta}\left(
X\right)  .$ \ We shall refer to these also as \emph{Koschorke classes}%
.\emph{\bigskip}

 There is an explicit formula expressing the $k_{p,q}$ in terms of the
classical Chern classes $c_{1},c_{2},...$ as Hankel determinants. \ The
proof is given in [K]. \ In particular, since we are interested in the
case where the index $k$ is zero, we have $p=q$, and $k_{p,p}$ has dimension
$2p^{2}.$ \ The first few Hankel determinants are:
\begin{align*}
h_{1,1}  & =c_{1}\\
h_{2,2}  & =\left\vert
\begin{array}
[c]{cc}%
c_{2} & c_{3}\\
c_{1} & c_{2}%
\end{array}
\right\vert =c_{2}^2-c_{3}c_{1}\\
h_{3,3}  & =\left\vert
\begin{array}
[c]{ccc}%
c_{3} & c_{4} & c_{5}\\
c_{2} & c_{3} & c_{4}\\
c_{1} & c_{2} & c_{3}%
\end{array}
\right\vert .
\end{align*}
In general $h_{p,q}$ is the determinant of the $q\times q$-matrix whose $(i,j)^{\rm th}$ entry is $c_{i-j+p}$.  \ (Here we interpret $c_0$ as 1, and $c_k$ as 0 when $k<0$.)

\bigskip

To see where these formulae come from, let us suppose that the formal series
\[c(t)=1+c_1t+c_2t^2+c_3t^3+ \ldots \]
is the expansion of a rational function
\[c(t)=\frac{a(t)}{b(t)}=\frac{1+a_1t+\ldots +a_pt^p}{1+b_1t+\ldots +b_qt^q}.\]
This means that each $c_n$ is expressed as a polynomial of degree $2n$ in the variables $a_1,\ldots ,a_p ;b_1,\ldots ,b_q$, when $a_k$ and $b_k$ are given weight $2k$. \ The resulting map of polynomial rings
\[\Q[c_1,c_2,c_3,\ldots ] \to \Q[a_1,\ldots ,a_p] \otimes \Q[b_1,\ldots ,b_q] \tag{9.1} \]
can be identified with the map of rational cohomology induced by the map
\[ BU_p \times BU_q \to {\rm Fred}_{p-q}(\H)\]
which classifies the formal difference of bundles $E_p \otimes 1 - 1\otimes E_q$, where $E_p$ and $E_q$ are the universal $p$- and $q$-dimensional bundles on $BU_p$ and $BU_q$, and $a_1,\ldots ,a_p$ and $b_1, \ldots ,b_q$ are the universal Chern classes in $H^*(BU_p)$ and $H^*(BU_q)$. We can identify $BU_n$ for any $n$ with the Grassmannian manifold of $n$-dimensional subspaces of $\H$, and then we have an obvious fibration
\[F_{p,q}\to BU_p \times BU_q, \tag{9.2}\]
whose fibre is the general linear group of Hilbert space, which is contractible.

It will also be helpful to introduce the embedding of polynomial rings
\[\Q[a_1,\ldots ,a_p] \otimes \Q[b_1,\ldots ,b_q] \to \Q[x_1,\ldots ,x_p] \otimes \Q[y_1,\ldots ,y_q] \tag{9.3}\]
defined by the factorizations
\begin{eqnarray*}
a(t) & = & (1+x_1t)\ldots (1+x_pt), \\
b(t) & = & (1+y_1t)\ldots (1+y_qt).
\end{eqnarray*}

The polynomial $R(a,b)$ of degree $2pq$ in the $a_i$ and $b_j$ that maps to $\prod (x_i-y_j)$ is called the {\it resultant} of the polynomials $a(t)$ and $b(t)$. It is well-known in classical algebra because its vanishing is the condition on the coefficients of the polynomials $a(t)$ and $b(t)$ for them to have a common factor.
\bigskip

The following proposition gives us the basic facts about Koschorke classes and Hankel determinants.
\bigskip

\noindent{\bf Proposition 9.4}

{\it (i) \ \ The homomorphism 9.1 is injective in degrees $<2(p+1)(q+1)$.

(ii) \ \ The Hankel determinant $h_{p,q}(c)$ maps to the resultant $R(a,b)$.

(iii) \ \ The determinant $h_{p',q'}(c)$ maps to zero if $p'>p$ and $q'>q$.

(iv) \ \ The Koschorke class $k_{p,q}$ is given, up to sign, by the Hankel determinant $h_{p,q}(c)$.}

\bigskip

{\bf Proof} \ \ We begin with (iii). \ The relation $c(t)b(t)=a(t)$ expresses the polynomial $a(t)$ as a complex linear combination of the $q+1$ formal power series 
$$ c(t),\ tc(t), \ t^2c(t), \ \ldots \ , t^qc(t).$$
In other words, the row vector
$$(1 \ \  a_1 \  \ a_2 \  \ \ldots \ a_p \  \ 0 \  \ 0 \  \ 0 \  \ \ldots \  )$$
is a linear combination of the $q+1$ rows
$$\begin{array}{ccccccc}
(1  \ \ & c_1 \ & c_2 \ & c_3 \ & \ldots \ & c_p \ & \ldots \ ) \\
(0 \ & 1 \ \ & c_1 \ & c_2 \ & \ldots \ & c_{p-1} \ & \ldots \ )\\
(0 \ & 0 \ \ & 1 \ \ & c_1 \ & \ldots \ & c_{p-2} \ & \ldots \ ) \\
  \cdots \\
(0 \ & 0 \ \ & 0\ \ & 0 \ \ & \ldots \ & c_{p-q} \ & \ldots \ ).
\end{array}$$

The Hankel determinant $h_{p',q+1}$ is the determinant of the $(q+1)\times (q+1)$-submatrix extracted from these rows, beginning with the entry $c_{p'}$ in the top left-hand corner.  \ But because $a_n=0$ if $n>p$ the rows of the submatrix are linearly dependent if $p'>p$, i.e. $h_{p',q+1}(c)=0$ if $p'>p$. \ The result (iii) follows, as  when $q'>q+1$ the determinant $h_{p',q'}  $ is that of a $q' \times q'$-submatrix of an array with $q'$ rows which are a fortiori linearly dependent.
\bigskip

To prove (i), consider the map (9.1) with $p$ and $q$ replaced by $p+k$ and $q+k$. \ Let $I_k$ be its kernel. \ Clearly $I_k$ is zero in degrees $\leq 2n$ if $p+k \geq n$. \ So if a polynomial $f(c)$ of degree $2n$ belongs to $I_0$ we can find $k\leq n-p$ such that $f(c)\in I_{k-1}$ but $f(c)\notin I_k$. \ Then $f(c)$ gives rise to a non-zero polynomial of degree $2n$ in $a_1,\ldots ,a_{p+k};b_1, \ldots ,b_{q+k}$ which vanishes if the polynomials $a(t)$ and $b(t)$ have a common factor. \ It is therefore divisible by $R(a,b)$, and so $n \geq (p+k)(q+k) \geq (p+1)(q+1)$, as we want.
\bigskip

To prove (ii), we consider the determinant $h_{p,q}(c)$ as a function of the coefficients of the rational function $a(t)/b(t)$.  \ If the polynomials $a(t)$ and $b(t)$ have a common factor, then the proof of (iii) just given, but with $q$ replaced by $q-1$, tells us that $h_{p,q}(c)$ will vanish. \ It must therefore be divisible by the resultant polynomial $R(a,b)$. \ As it has the same degree $2pq$ as $R(a,b)$, it must be a rational multiple of it. \ But if we specialize by setting $a_i=0$ and $b_j=0$ except for $a_p$, then both $R(a,b)$ and $h_{p,q}$ become $(a_p)^q$.
\bigskip

Finally, to prove (iv), it is enough by (i) to show that $k_{p,q}$ and $h_{p,q}(c)$ have the same restriction to $F_{p,q}$. \ But let us recall that the cohomology class represented by $\bar{F}_{p,q}$ becomes, when restricted to $F_{p,q}$, the Euler class of the normal bundle to $F_{p,q}$. \ This normal bundle is $(E_p \otimes 1)^*\otimes (1 \otimes E_q)$, and its Euler class is $\prod (y_j -x_i)$. $\Box$

\bigskip

The geometry shows that the Koschorke classes extend from the fibre Fred$(\H)$ to the
total space $Y$ of the universal bundle, and hence that they are \emph{invariant} in the sense of Section 8. \  But we can also see the invariance from Proposition 9.4. \ The maps (9.1) and (9.3) are equivariant with respect to the action of the 1-parameter group of section 8, which acts on the variables $x_i$ and $y_j$ of (9.3) by $x_i \mapsto x_i +u$ and $y_j \mapsto y_j + u$. \ The resultant $\prod (x_i - y_j)$ is clearly invariant under  this action, and so, therefore, is the class $h_{p,q}(c)$. \ It is worth remarking that the action of the 1-parameter group on the total Chern class $c(t)$ of a class of virtual dimension 0 is $c(t) \mapsto c(\tilde{t})$, where $\tilde{t}=t/(1-ut)$.  
\bigskip

It is clearly attractive to have such nicely defined geometric invariant
classes, but unfortunately they are far from generating all invariant classes.
\ They are too thin on the ground, as we see by noticing that their dimensions
$2p^{2}$, i.e.  \;  2,\;8,\;18,\;... \;   rise rapidly. \  In fact they are invariant not only under the 1-parameter group described in (8.6) but even under an action of the infinite dimensional group of formal reparametrizations
\[t\mapsto \tilde{t}=t + u_1t^2 +u_2t^3+u_3t^4 + \cdots \; , \]
of the line, in which the previous group sits as the subgroup of M\"{o}bius transformations $$t\mapsto t/(1-ut).$$
This is proved in a computational way in [Se2] Prop.6.4, and unfortunately we do not know any conceptual proof. \ One complicating thing is that the larger reparametrization group does not act on the cohomology of Fred$(\H)$ by ring automorphisms, and so it does not change the Chern classes simply by $c(t)\mapsto c(\tilde{t})$, but by a twisted version of that.

\bigskip 

 It is not clear what special role the Koschorke classes play, but we believe that they deserve further study.

\bigskip

\section{Operations in twisted K-theory}

 \ \ \ \ Ordinary $K$-theory has operations arising from representation
theory. \ Thus we have the exterior powers $\lambda^{k},$ the symmetric powers
$s^{k}$, and the Adams operations $\psi^{k}.$ \ In [A1] the ring of
operations is described in terms of the decomposition of the tensor powers of a bundle 
under the symmetric groups. \ If we define
\[
R_{\ast}=\sum\limits_{k}\text{Hom}_{\Z}\left(  R\left(  S_{k}\right)
,\Z\right)
\]
where $R\left(  G\right)  $ is the character ring of a group $G$, then there is a
natural homomorphism
\[
R_{\ast} \to {\rm Op}\left(  K\right)
\]
which is an embedding . \ If we think of ${\rm Op}\left(  K\right)  $
in terms of self-maps of the classifying space Fred$\left(  \H\right)  $ of $K$
then it acquires a topology from the skeleton filtration, and it is not hard to
see that $R_{\ast}$ has dense image in ${\rm Op}\left(  K\right)  ,$ so that
${\rm Op}\left(  K\right)  $ can be viewed as the completion $\hat{R}_{\ast}%
.$\emph{\bigskip}

 In [A1] it is shown that $R_{\ast}$ is a polynomial ring on the
$\lambda^{k},$ or on the $s^{k}.$ \ Over $\Q$, the $\psi^{k}$ are also  polynomial generators.
\emph{\bigskip}

 There is actually another natural set of operations not included in
$R_{\ast}$. \ These involve the operation $\ast$ given by taking the dual
$E^{\ast}$ of a vector bundle $E.$ \ It is natural to denote $\ast$ by
$\psi^{-1},$ as can be seen by computing the Chern character of $\psi^{k}$
applied to the universal bundle $\xi$:

\noindent if
\[
{\rm ch}\,\xi=\sum\limits_{n}{\rm ch}_{n}\,\xi
\]
then
\[
{\rm ch}\, \psi^{k}\xi=\sum\limits_{n}k^{n}{\rm ch}_{n}\,\xi
\]
while
\[
{\rm ch}\, \xi^{\ast}=\sum\limits_{n}\left(  -1\right)  ^{n}{\rm ch}_{n}\,\xi
\]
The usual composition rules
\[
\psi^{k}\ {\small o}{\tiny \ \ }\psi^{l}=\psi^{kl}=\psi^{l}\ {\small o}%
{\tiny \ \ }\psi^{k}%
\]
for positive integers $k,l$ then extend also to negative integers when we
define (for $k>0$)
\[
\psi^{-k}=\psi^{k}\ {\small o\ \ }\psi^{-1},%
\]
i.e.
\[
\psi^{-k}\left(  E\right)  =\psi^{k}\left(  E^{\ast}\right)  =\left[  \psi
^{k}\left(  E\right)  \right]  ^{\ast}.
\]
\textbf{Note \ }We must distinguish carefully between the \emph{composition}
$\psi^{k}\ {\small o\ \ }\psi^{l}$ and the \emph{product} $\psi^{k}\psi^{l}.$
\ Thus%
\[
\psi^{k}\ {\small o\ \ }\psi^{l}\left(  E\right)  =\psi^{k}\left(  \psi
^{l}E\right)  =\psi^{kl}E
\]
is quite different from the product $\psi^{k}\left(  E\right)  \psi
^{l}\left(  E\right)  $.

\bigskip

 The operation $\psi^{-1}$ automatically belongs to the completed ring $\hat{R}_{\ast}$, and in fact 
\[\psi^{-1}(x) = - \left(\sum\gamma^p(x)\right)^{-1}\sum p\gamma^p(x)\]
in terms of the operations $\gamma^p$ defined below.
Thus if we allow ourselves to work within the completion
$\hat{R}_{\ast}$  there is nothing new to be gained from the $\psi^{k}$
for negative $k,$ but if we want to avoid completion then the negative
exponents do give something new. \ This will become important shortly when we
pass to the twisted $K$-theory.\bigskip

 The treatment in [A1] is in terms of complexes of vector bundles,
but it extends naturally to Fredholm complexes [S1]. \ Here a single
Fredholm operator
\[
T:\H\rightarrow \H
\]
is generalized to a complex of Hilbert spaces
\[
0\rightarrow \H_{0}\overset{T_{0}}{\rightarrow}\H_{1}\overset{T_{0}}%
{\rightarrow}\H_{2}...\rightarrow \H_{n}\rightarrow0
\]
whose homology is finite-dimensional. \ We can still take tensor products and
decompose under the symmetric group.\bigskip

 Now we are ready to consider the twisted case. \ A twisted family of
Fredholm operators over a base space $X$ can be tensored with itself $n$ times
to give a twisted family of Fredholm complexes over $X$, and then decomposed
under the symmetric group $S_{n}.$ \ The important point to note is that the
twisting element $\zeta$ we get is now $n$ times the original twisting element
$\eta.$ \ This  follows from Proposition (2.1)(vii) in [AS].
\ Thus it would seem that (for $\eta$ of infinite order) we always get operations
that shift the twisting, and so are not \textquotedblleft
internal\textquotedblright. \ This is where the duality operation $\ast$ or
$\psi^{-1}$ comes to our rescue, as we shall now explain.\bigskip

 Consider any monomial $\mu$ in the $\psi^{k}$:
\[
\mu=\psi^{k_{1}}\psi^{k_{2}}...\psi^{k_{l}}%
\]
where the total degree $\sum k_{i}=1.$ \ Of course this requires that some of
the $k_{i}$ are \emph{negative}. \ We can apply this to our Fredholm family
with twisting element $\eta$ and get a new Fredholm family with twisting
element
\[
\left(  \sum k_{i}\right)  \eta=\eta ,
\]
as follows by using Proposition (2.1)(viii) of [AS]. \ Thus the
monomial $\mu$ does define an \emph{internal} operation which preserves the
twisting.\bigskip

 It seems reasonable to expect that in this way we get (over $\Q$) a
dense set of all internal operations. \ However, it is not obvious how to
establish this. \ The basic reason is that the filtration that is normally
used in $K$-theory is generated by the Grothendieck $\gamma^{k}$-operations
which are defined in terms of $\lambda_{t}=\sum\lambda^{i}t^{i}$ by
\[
\gamma_{t}=\sum\gamma^{k}t^{k}=\lambda_{t/1-t}.
\]
This means that $\gamma^{n}$ is a linear combination of $\lambda^{k}$ for
$k\leq n,$ and so in the twisted case it mixes up the twistings. \ We cannot
therefore use the $\gamma^{n}$ to give a filtration on our internal operations
for twisted $K$-theory. \ This might of course suggest that our whole question
is unnatural and that one should consider only the $K$-theory got by summing
over all twistings.\bigskip

 It also seems clear that we should not need all the monomials $\mu$
of total degree $1$ - they would be over-abundant. \ One might, for example,
hope to restrict oneself to those $\mu$ with at most one negative exponent,
say $k_{l}$:\   these are  determined by the arbitrary string of positive exponents
$k_{1},k_{2},...k_{l-1}.$ \ We would also include the case
$\mu=\psi^{1}$ (the identity), where there is no negative exponent.\bigskip

 Since we are only interested in elements of augmentation or index
$0,$ a monomial $\mu$ will raise the filtration from $1$ to $l.$\bigskip

 We now revert to the question raised in Section 8 of finding
natural lifts of the invariant cohomology classes into twisted cohomology
classes. \ We now have the twisted classes given by the $\mu$ and taking the
leading term (i.e. dimension $2l$) of their Chern characters we must end up
with invariant classes. \ Thus our internal operations, followed by the Chern
character, give in principle some natural lifts. \ In other words we have a
method to construct characteristic classes which end up in the twisted
cohomology. \ It remains a purely algebraic problem to find some natural
basis.\bigskip

 Note that we could use monomials in the $\lambda^{k}$%
\[
\lambda^{k_{1}}\lambda^{k_{2}}...\lambda^{k_{l}},%
\]
with $\sum k_{i}=1$, in exactly the same way, if we define $\lambda^{k}$
for negative $k$ by
\[
\lambda^{-k}\left(  \xi\right)  =\lambda^{k}\left(  \xi^{\ast}\right)  =
\left[  \lambda^{k}\left(  \xi\right)  \right]  ^{\ast}.
\]
These are more likely to produce some form of integral basis, but they are much
more difficult to handle algebraically because their Chern characters are
complicated.\bigskip

 Finally we make a remark about the Koschorke classes. \ Since these
are defined by complex submanifolds $F_{p,q}$ of Fred$\left(  H\right)  $,
their closures $\bar{F}_{p,q}$ should define complex analytic subspaces and
hence, by using a resolution of the ideal-sheaf, an element of $K\left(
\text{Fred}\left(  \H\right)  \right)  .$ \  But this requires some care, since
we are in infinite dimensions, and the approximations are more difficult to
handle in $K$-theory because it is only filtered and not graded. \ However,
assuming these difficulties are overcome, we shall end up with classes
$\tilde{k}_{p,q}$ in the ordinary (untwisted) $K$-theory of the base. \ These
provide \textquotedblleft lifts\textquotedblright\ for the Koschorke classes
in the Atiyah-Hirzebruch spectral sequence.

\section{Appendix}

\ \ \ \ In this appendix we shall give simple examples of spaces with
non-zero Massey products of the kind occurring in Section 6 as higher
differentials in the Atiyah-Hirzebruch spectral sequence of twisted
$K$-theory. \ We are indebted to Elmer Rees for these examples. \ We shall
also discuss some general issues arising from them.

\bigskip

 We begin by introducing a 3-manifold $Y$ which has a non-zero Massey
triple product $\left\{  x,x,y\right\}  $ with $x,y\in H^{1}\left(
Y\right)  .$ \ We will then take the product $Y\times \C P^{2}$ with the complex
projective plane. \ This will have the Massey triple product $\left\{
xt,\ xt,\ y\right\}  $ non-zero, where $t$ generates $H^{2}\left(
\C P^{2}\right)  $.

\bigskip

 To contruct $Y$ we just take the $U(1)$-bundle over a 2-torus
$T=S^{1}\times S^{1}$ with Chern class 1. \ If $x,y$ are closed 1-forms on $T,
$ coming from the two factors, and if $z$ is the standard connection 1-form on
$Y$ (normalized so that its integral over the $U\left(  1\right)  $-fibre is
1), then $dz=xy$ gives the (normalized) curvature 2-form. \ The cohomology of
$Y$ has the following generators.%
\begin{equation}
\left.
\begin{array}
[c]{lll}%
1 & \text{in} & H^{0}\\
x,y & \text{in} & H^{1}\\
xz,yz & \text{in} & H^{2}\\
xyz & \text{in} & H^{3}%
\end{array}
\right\} \tag{A1}%
\end{equation}
We see that the Massey triple product%
\[
\left\{  x,x,y\right\}  =xz\ \ \in H^{2}%
\]
is non-zero. \ Moreover, it is not in the image of multiplication by $x$ (note
that $z$ is \textbf{not} a closed form). \ In fact multiplication by $x $
annihilates all of $H^{1}.$

\bigskip

 Now consider $M=Y\times \C P^{2}.$ \ The generators in $(A1)$ get
multiplied by $1,t,t^{2}$ to give the generators of $H^{\ast}\left(  M\right)
$. \ We see that%
\begin{equation}
\left\{  xt,xt,y\right\}  =xzt^{2}\ \in H^{6}\tag{A2}%
\end{equation}
Moreover this class is not in the image of multiplication by $xt$ from
$H^{3}$ to $H^{6}.$ \ In fact $xt$ annihilates all three generators
$xt,yt,xyz$ of $H^{3}.$

\bigskip

 Thus on $M,$ for the twisting class $xt,$ the differential 
\[
d_{5}:E_{4}\longrightarrow E_{4}%
\]
of the Atiyah-Hirzebruch spectral sequence is non-zero on the element $y\in E_{4}.$

\bigskip

 A little computation shows in fact that the ranks of $E_{2}%
,\ E_{4},\ E_{6}=E_{\infty}$ are 18, 14, 10 respectively. \ Thus the rank of
$K_{P}^{\ast}\left(  M\right)  $ is 10 (where $\left[  xt\right]  $ is the
class of $P$). \ This can also be checked directly by computing the cohomology
of the operator $d+xt$ on differential forms as in Section 6.

\bigskip

 Note that $Y$ is just the Eilenberg-MacLane space $K(\pi ,1)$ for the Heisenberg
group $\pi$ generated by three elements $a,b,c$ with $c$ central and equal to the
commutator of $a,b.$ \ The group $\pi$ is a central extension of $\Z^{2}$ by $\Z.$ \ All the
computations above are just computations in the cohomology of this discrete group.

\bigskip

 We can generalize this example to a tower of examples, each a
$U\left(  1\right)  $-bundle over its predecessor with Chern class the
previous Massey product. \ The first step beyond $Y$ gives a 4-manifold $Z$
with Chern class $xz,$ having a connection 1-form with $xz=du.$ \ The Massey
product%
\[
\left\{  x,x,x,y\right\}  =xu
\]
is then non-zero. \ The construction provides a sequence of examples in which
the successive Massey products are, at each stage, non-zero. \  Note that all these examples arise from discrete groups which are
successive central extensions by $\Z,$ so that again all calculations are just
calculations in the cohomology of discrete groups.

\bigskip

 A more direct construction of the $(n+1)$-dimensional manifold $Y_n$ in this tower is as a bundle over the circle whose fibre is the $n$-dimensional torus $T^n = \R^n/\Z^n$, and whose monodromy is the integral unipotent matrix $1+N$, where $N$ is the Jordan block with 1s in the superdiagonal. Then $H^1(Y_n)=\Z^2$, generated by $x$ coming from the base circle, and $y$ which restricts to the generator of $H^1(T^n)$ invariant under the monodromy. \ In $Y_n$ the $(n+1)$-fold Massey product $\{x,\ldots,x,y\}\in H^n(Y_n)$ is non-zero.
\bigskip

 To get to Massey products involving a class in $H^{3}$ we again
multiply by a complex projective space of the relevant dimension. \ Thus for
$Y_3$ we use $Y_3\times \C P^{3}$ and find%
\[
\left\{  xt,xt,xt,y\right\}  =xut^{3}\in H^{8},%
\]
and verify that  the  element 
\[
d_7\left(  y\right)  =xut^{3}%
\]
in the $E_6$-term of the spectral sequence is non-zero.

\bigskip

 Clearly, in all these examples, we could replace the finite dimensional $\C P^{n}
$ by $\C P^{\infty}$ and we could interpret our results in terms of equivariant
calculations for the group $G=U\left(  1\right)  .$ \ The class $xt $ that we
have been using is just the equivariant class generating $H_{G}^{3}(G)$ that
was studied in Section 6 of [AS] (note that $G$ acts on $G$ by conjugation, and so trivially when $G=U\left(  1\right)  $). \ We observed there that this
is the class of a natural $P\left(  \H\right)  $-bundle over $G=U\left(
1\right)  ,$ where $\H$ is a graded Hilbert space, and the grading shifts by 1 on flowing round the circle. \ Our example therefore deals with
this equivariant bundle $P$ on $S^{1}$ lifted from the first factor of $T^{2}$
to $Y.$ \ The example is therefore not as \emph{ad hoc} as it might seem.

\bigskip

 Instead of treating $t$ as an equivariant cohomology class we can
also  treat it simply as a real or complex variable. \ Then the operator%
\[
D_{t}=d+xt
\]
can be viewed as a covariant derivative for a connection, depending on a
parameter $t.$ \ Regarding $t$ as constant, we have $D_{t}^{2}=0$, so that the
connection is flat. \ Its cohomology is then just the standard cohomology of a
flat connection. \ When we regard $t$ as a variable the cohomology of $D_{t}$
becomes naturally a module over the polynomial ring $\C[t].$

\bigskip

 Since $S^{1}\times \C P^{\infty}$ is the universal space for pairs of
integer cohomology classes $x\in H^{1},$ $t\in H^{2},$ the interpretation of
twisted cohomology in terms of the cohomology of flat connections applies
whenever our twisting class $\eta$ decomposes as $xt.$ \ More generally, the
same applies whenever $\eta$ is a linear combination $\eta=\sum x_{i}t_{i}.$
\ In particular, this occurs in the following situation related to the twisted
equivariant $K$-theory $K_{G,P}\left(  G\right)  $ in Section 6 of [AS],
where $G$ is a compact Lie group acting on itself by conjugation. \ The  case we have just been examining is when $G=U(1)$, so let us now
restrict to the case where $G$ is simply-connected. \ In [A2] the ordinary
$K$-theory of $G$ was usefully studied via the natural map%
\[
\pi:G/T\ \times\ T\longrightarrow G,
\]
where $T\ $ is a maximal torus, and $\pi\left(  gT,u\right)  =gug^{-1}$.
\ Clearly this map is $G$-equivariant for the conjugation action on $G,$ so it
induces a map in $G$-equivariant $K$-theory or $G$-equivariant
cohomology. \ Since%
\[
H_{G}^{\ast}\left(  G/T\right)  =H_{T}^{\ast}\text{(point)}=H^{\ast}\left(
BT\right),
\]
we get a homomorphism%
\[
H^{\ast}(B_{T})\otimes H^{\ast}\left(  T\right)  \overset{\pi^{\ast}%
}{\longleftarrow}H_{G}^{\ast}(G).
\]
The twisting class%
\[
\left[  \eta\right]  \in H_{G}^{3}\left(  G\right)
\]
that we introduced in [AS] then necessarily gives a decomposable element%
\[
\pi^{\ast}[\eta]=\sum x_{i}t_{i},%
\]
where $t_{i}\in H^{1}\left(  T\right)  $ and $x_{i}\in H^{2}\left(  BT\right).
$ \ (It is easy to see there can be no term in $H^{3}(T)$). \ Thus, the twisted
equivariant $K$-theory of $G/T\times T$ (twisted by $\pi^{\ast}\eta)$ can be
understood in terms of families of flat connections, parametrized by $\left(
t_{1},...,t_{n}\right)  .$

\bigskip

 A crucial point, exploited in [A2], is that $\pi$ preserves
orientation and dimension. \ (Its degree is the order of the Weyl group
of $G$). \ This implies, purely formally, that $\pi^{\ast}$ on real cohomology
is injective and canonically split. \ The same argument applies to twisted
(real) cohomology, and to twisted $K$-theory (tensored with $\R$). \ It also
extends to the equivariant case. \ Thus we can view the twisted equivariant
$K$-theory of $G$ (modulo torsion) as a subspace of the corresponding
group for $G/T\times T.$ \ The corresponding cohomology is therefore
interpretable in terms of families of flat connections.

\bigskip

 The groups $K_{G,P}(G)$ are the subject of the work by Freed,
Hopkins, and Teleman [FHT], referred to in [AS], and are shown by them to be related
to the Verlinde algebra. \ \ This in turn is connected to conformal field
theory and to Chern-Simons theory. \ There were speculations in [A3] that
Chern-Simons theory for $G$ might be \textquotedblleft
abelianized\textquotedblright, i.e. reduced to its maximal term $T.$ \ This
programme has been carried out in detail by Yoshida [Y]. \ There might be
some connection between this programme and the process described above
replacing $G$ by $G/T\times T.$

\section*{References}

\begin{enumerate}

\item[[A1]]\ \ M.F.Atiyah,\   Power operations in $K$-theory.  \ {\it Quarterly J. Math., Oxford} \ (2), \ {\bf 17}  \ (1966), \  165 -- 193.
\item[[A2]]\ \ M.F.Atiyah, \ On the $K$-theory of compact Lie groups.  \ {\it Topology} \  {\bf 4}  \ (1965), \  95 --99.
\item[[A3]] \ \ M.F.Atiyah, \ {\it The Geometry and Physics of Knots.} \ Lincei Lectures, Cambridge Univ. Press, 1990. 
\item[[AH]] \ \ M.F.Atiyah and F.Hirzebruch, \ Vector bundles and homogeneous spaces. \  {\it Proc. Symp. Pure Mathematics}  \ {\bf 3}, Amer. Math. Soc. 1961.
\item[[AS]]\ \ M.F.Atiyah and G.B. Segal, \  Twisted $K$-theory. \ {Ukrainian Math. Bull.} \  {\bf 1} \ (2004).
\item[[BCMMS]] \  P. Bouwknegt, A.L.Carey, V. Mathai, M.K.Murray, and D. Stevenson, \  Twisted $K$-theory and $K$-theory of bundle gerbes. \  {\it Comm. Math. Phys.}\ {\bf 228} \  (2002), 17--45.
\item[[BGV]]\ \ N. Berline, E.Getzler, and M.Vergne,  \ {\it Heat Kernels and Dirac operators.} \  Springer, Berlin, 1992.
\item[[DGMS]]\ \ P. Deligne, P. Griffiths, J. Morgan, and D. Sullivan, \ Real homotopy theory of K\"ahler manifolds. \ {\it Invent. Math.}, \ {\bf 29} \  (1975), 245 -- 274.
\item[[FHT]]\ \ D. Freed, M. Hopkins, and C. Teleman, \ Twisted equivariant $K$-theory with complex coefficients. \ math.AT/0206257.
\item[[K]]\ \ U. Koschorke, \ Infinite dimensional $K$-theory and characteristic classes of Fredholm bundle maps. \ {\it Proc. Symp. Pure Math.} \  {\bf 15} \  Amer. Math. Soc. 1970, 95 --133.
\item[[Q]]\ \ D.G.Quillen, \ Superconnections and the Chern character.\ {\it Topology}  \ {\bf 24} \  (1985), 89 -- 95.
\item[[R]]\ \ J. Rosenberg, \ Continuous-trace algebras from the bundle-theoretic point of view.  \ {\it J. Austral. Math. Soc.} \  {\bf A47} \  (1989), 368--381.
\item[[S1]]\ \ G.B.Segal, \ Fredholm complexes.\ {\it Quarterly J. Math., Oxford}   {\bf 21} \  (1970), 385--402. 
\item[[S2]]\ \ G.B.Segal, \  Unitary representations of some infinite dimensional groups.\ {\it Comm. Math. Phys.} \ {\bf 80} \  (1981), 301--342.
\item[[S3]] \ \ G.B.Segal, \ Topological structures in string theory.  \ {\it Phil. Trans. Roy. Soc. London}   \  {\bf A359} (2001), 1389 -- 1398.
\item[[Su]]\ \ D. Sullivan, \ Infinitesimal computations in topology. \ {\it Inst. Hautes \'Etudes Sci. Publ. Math.} \  {\bf 47} \  (1977), 269 -- 331.
\item[[TX]] \ \ J.-L. Tu and P. Xu, \ Chern character for twisted $K$-theory of orbifolds. \ math.KT/0505267. 
\item[[Y]]\ \ T. Yoshida,  \ An abelianization of SU(2) WZW model. \ {\it Annals of Math.} \  (2005) \ (to appear).
\end{enumerate}

\end{document}